\documentclass[12pt]{amsart}
\usepackage{amssymb}
\usepackage[cmtip,matrix,arrow]{xy}

\newtheorem{theorem}{Theorem}[section]

\newtheorem{proposition}[theorem]{Proposition}

\newtheorem{lemma}[theorem]{Lemma}
\theoremstyle{definition}    
\newtheorem{definition}[theorem]{Definition}
\theoremstyle{remark}

\newtheorem{remark}[theorem]{Remark}

\newtheorem{example}[theorem]{Example}
\newcommand{\labell}\label

\newcommand{\beq}{\begin{eqnarray*}}
\newcommand{\eeq}{\end{eqnarray*}}
\renewcommand{\H}{{\mathcal H}}

\newcommand{\W}{\mathcal{W}}

\renewcommand{\L}{\mathcal{L}}
\newcommand{\La}{\Lambda}
\renewcommand{\O}{\mathcal{O}}
\newcommand{\T}{\mathcal{T}}
\newcommand{\J}{\mathcal{J}}
\newcommand{\Co}{\mathcal{C}}
\newcommand{\ca}{\mathcal}

\renewcommand{\c}{\mathfrak{c}}

\newcommand{\F}{\mathcal{F}}
\newcommand{\N}{\mathbb{N}}
\newcommand{\R}{\mathbb{R}}
\newcommand{\C}{\mathbb{C}}

\newcommand{\SU}{\on{SU}}

\newcommand{\Z}{\mathbb{Z}}

\DeclareMathOperator{\Cl}{Cl}
\newcommand\lie[1]{\mathfrak{#1}}

\newcommand{\g}{\lie{g}}
\newcommand{\m}{\lie{m}}

\renewcommand{\t}{\lie{t}}

\newcommand{\on}{\operatorname}

\newcommand{\Eul}{ \on{Eul} } 
 
\newcommand{\Ad}{ \on{Ad} }
\newcommand{\ad}{\on{ad}}

\renewcommand{\ker}{ \on{ker}}

\newcommand{\Spin}{ \on{Spin}}
\newcommand{\SO}{ \on{SO}}
\newcommand{\Mult}{  \on{Mult}}
\newcommand{\Vol}{  \on{Vol}}

\newcommand{\diag}{  \on{diag}}

\newcommand{\tpi}{{2\pi i}}

\newcommand\qu{/\kern-.7ex/} 

\newcommand{\hra}{\hookrightarrow}

\renewcommand{\d}{{\mbox{d}}}
\newcommand{\ol}{\overline}

\newcommand\Phinv{\Phi^{-1}}

\newcommand\sig{\sigma}

\newcommand\Om{\Omega}
\newcommand\om{\omega}

\newcommand{\f}{\frac}

\newcommand{\olt}{\overline{\theta}}
\newcommand{\p}{\partial}
\renewcommand{\l}{\langle}
\renewcommand{\r}{\rangle}
\newcommand{\bfl}{{\bf \,\langle\!\langle\,}}
\newcommand{\bfr}{{\bf \,\rangle\!\rangle\,}}
\newcommand\hh{{\f{1}{2}}}

\newcommand{\pt}{\on{pt}}

\newcommand{\wh}{\widehat}

\newcommand{\mf}{\mathfrak}
\renewcommand{\S}{\Delta}   
\begin{document}

\title[Localization]{Group-valued equivariant localization}
\author{A. Alekseev}
\address{Institute for Theoretical Physics \\ Uppsala University \\
Box 803 \\ \mbox{S-75108} Uppsala \\ Sweden}
\email{alekseev@teorfys.uu.se}

\author{E. Meinrenken}
\address{University of Toronto, Department of Mathematics,
100 St George Street, Toronto, Ontario M5R3G3, Canada }
\email{mein@math.toronto.edu}

\author{C.Woodward}
\address{Mathematics -- Hill Center, Rutgers University,
110 Frelinghuysen Road, Piscataway NJ 08854-8019, USA}
\email{ctw@math.rutgers.edu}

\date{May 1999}

\begin{abstract}
We prove a localization formula for group-valued equivariant de Rham
cohomology of a compact $G$-manifold. This formula is a non-trivial
generalization of the localization formula of Berline-Vergne and
Atiyah-Bott for the usual equivariant de Rham cohomology. As an
application, we obtain a version of the Duistermaat-Heckman formula
for Hamiltonian spaces with group-valued moment maps.
\end{abstract}

\maketitle

\section{Introduction}
\label{Sec:Intro}
The main result of this paper is a localization theorem 
for $G$-manifolds $M$ which we view as the group valued
analogue of the localization formula (Berline-Vergne
\cite{be:ze}, Atiyah-Bott \cite{at:mom}) in equivariant de Rham theory. 

Recall that the equivariant cohomology of $M$ (over $\C$)
is the cohomology of 
Cartan's \cite{ca:no} complex $C_G(M)=(S\g_\C^*\otimes  \Om(M))^G$
of equivariant differential forms, equipped with 
a certain differential $\d_G$. The complex $C_G(M)$ embeds 
into the larger complex complex 
$\wh{C}_G(M)=(\ca{E}'(\g^*)\otimes \Om(M))^G$ where 
$\ca{E}'(\g^*)$ is the convolution algebra of compactly supported 
distributions, and $S\g^*$ is viewed as distributions 
supported at $0$. If $M$ is compact and oriented, 
integration over $M$ defines a chain map 
$$ \int_M:\ \wh{C}_G(M)\to \wh{C}_G(\pt) = \ca{E}'(\g^*)^G.$$
For any cocycle $\beta\in \wh{C}_G(M)$, the localization formula
expresses the value of the Fourier transform 
$\l\int_M\beta,e^{\tpi\l\cdot,\xi\r}\r$ 
in terms of integrals over the connected 
components $F$ of the zero set of the vector field $\xi_M$
generated by $\xi$:
\begin{equation}\label{localization}
\l\int_M\beta,e^{\tpi\l\cdot,\xi\r}\r=
\sum_{F\in\ca{F}(\xi)} \int_F 
\f{\iota_F^* \l \beta,e^{\tpi\l\cdot,\xi\r}\r}{\Eul(\nu_F,\tpi\xi)}.
\end{equation}
Here $\Eul(\nu_F,\tpi\xi)$ denotes the equivariant Euler form of the
normal bundle. A special case of the localization theorem is the 
Duistermaat-Heckman ``exact stationary phase'' formula
from symplectic geometry, in which $\beta$ is the equivariant 
Liouville form. 

In \cite{al:no} we introduced a new version of the 
equivariant de Rham complex where the dual of the Lie algebra 
$\g^*$ is replaced by the Lie group $G$. That is, the chain complex 
is the space $\wh{\ca{C}}_G(M)=(\ca{E}'(G)\otimes\Om(M))^G$, 
with a suitable differential $\d_G$ depending on the choice
of an invariant inner product on $\g$. 
Integration is a chain map
$$ \int_M:\,\wh{\ca{C}}_G(M)\to \wh{\ca{C}}_G(\pt)=\ca{E}'(G)^G.$$

Our localization formula describes the Fourier coefficients 
of the integral $\int_M\beta$ of any equivariant cocycle,
that is, its pairings with irreducible characters of $G$. 
Let $T$ be a maximal torus of $G$, with Lie algebra $\t$, and let 
$\lambda\in \t^*$ be a dominant weight for some choice of 
positive Weyl chamber. It parametrizes an irreducible 
representation $V_\lambda$, with character $\chi_\lambda$. 
Let $\rho\in\t^*$ be the 
half-sum of positive roots and use the inner product to identify
$\t^*\cong\t$. 
Our main result expresses the 
Fourier coefficient $\l \int_M\beta,\chi_\lambda\r $ in terms of 
integrals over connected components $F$ of the zero set 
of the vector field $(\lambda+\rho)_M$, as follows. 
\begin{equation}\label{q-localization} 
\l \int_M\beta,\chi_\lambda\r =\dim V_\lambda 
\sum_{F\in\ca{F}(\lambda+\rho)} 
\int_F 
\f{\iota_F^*\big( e^{\hh\iota({r}_M)}
\l\beta,\S_\lambda\r\big)}{\Eul(\nu_F,\tpi(\lambda+\rho))}.  
\end{equation}
It involves the bi-vector field $r_M$ defined by the classical $r$-matrix, 
and certain  ``spherical harmonics'' $\S_\lambda\in C^\infty(G)$ 
determined by $\lambda$. 

There are several conceptual differences between the localization
formulas \eqref{localization} and \eqref{q-localization}. In
\eqref{localization} the equivariant parameter $\xi$ varies
continuously, and for an open dense set of $\xi\in\t$ the fixed point
set is just the fixed point set for the maximal torus. It follows by
continuity that if the maximal torus has no fixed points then the
integral of an equivariant cocycle vanishes.  The situation in
\eqref{q-localization} is different because $\lambda$ only varies in a
discrete set. Hence the continuity argument does not apply, and
integrals of equivariant cocycles may be non-zero even if the maximal
torus has no fixed points.  Next, for $\xi=0$ the localization formula
\eqref{localization} becomes tautological, whereas
\eqref{q-localization} gives a non-trivial result even for the
$\lambda=0$ Fourier coefficient, which localizes to the set of zeroes
of $\rho_M$.

We use the new localization theorem to obtain Duistermaat-Heckman type
formulas for the theory of group-valued moment maps
\cite{al:mom,al:du}. In this case the localization contributions 
simplify and lead to formulas for intersection pairings on reduced
spaces. 

This result applies in particular to intersection pairings on 
moduli spaces of flat $G$-connections on surfaces. 
According to \cite{al:mom,al:du} these spaces are reductions of group 
valued Hamiltonian $G$-spaces. 
In this example, our localization formula produces the formulas 
conjectured by Witten in \cite{wi:tw} and confirmed (for $G=\SU(n)$) 
by Jeffrey-Kirwan \cite{je:in}. Details of the application to 
moduli spaces will be presented in a forthcoming paper.

The paper is organized as follows. Section \ref{Sec:EqCoh} is a review
of equivariant de Rham theory. In particular, we describe the ``group
valued'' version introduced in \cite{al:no}.  In Section
\ref{Sec:Abel} we construct a suitable ``restriction map'' from
$G$-valued equivariant cohomology to $T$-valued equivariant
cohomology, and in Section \ref{Sec:PropRes} we discuss some of its
properties.  This ``abelianization'' procedure is the key step in our
proof of the $G$-valued localization formula in Section
\ref{Sec:LocForm}. In Section \ref{Sec:Appl} we use the localization
formula to prove a Duistermaat-Heckman principle for 
group valued moment maps.

\section{Equivariant de Rham theory}
\label{Sec:EqCoh}

In this Section we recall the Cartan and Weil models of equivariant de
Rham theory, and their non-commutative versions introduced in
\cite{al:no}.  For a more detailed account of equivariant de Rham
theory see e.g. the forthcoming book of Guillemin-Sternberg
\cite{gu:su}.

Throughout this paper, $G$ will denote a compact, connected Lie 
group with Lie algebra $\g$. (Occasionally we will make some 
additional assumptions on $G$.)
 
For any manifold $M$ we denote by $\Om^\star(M)$ the de Rham complex
of complex-valued differential forms.  Given a $G$-action on $M$, let
$$\xi_M=\f{\d}{\d t}\big|_{t=0}\exp(-t\xi)^*$$ 
denote the fundamental vector field corresponding to $\xi\in\g$.  Let
$$L_\xi:\,\Om^\star(M)\to \Om^\star(M), \ \ \ 
\iota_\xi:\,\Om^\star(M)\to \Om^{\star-1}(M)$$ 
be the Lie derivative (resp. contraction) by $\xi_M$.  Given a basis
$\{e_a \}\subset \g$ we also use the notation $L_a=L_{e_a}$ and
$\iota_a=\iota_{e_a}$.

\subsection{Cartan models}
\subsubsection{Equivariant cohomology} 
Let $\ca{E}'(\g^*)$ be the convolution algebra of compactly supported
complex-valued distributions on $\g^*$. The subalgebra of
distributions supported at $0$ is identified with the symmetric
algebra $S\g^*_\C$, which is generated by the elements
$$ v^a= \frac{\d}{\d t} \ \delta_{t e^a} \Big|_{t=0}
= -\frac{\p}{\p \mu_a} \ \delta_0$$
corresponding to the dual basis vectors $e^a$.  Let
$$ C_G(M)=(S\g^*_\C\otimes \Om(M))^G,\ \
\wh{C}_G(M)=(\ca{E}'(\g^*)\otimes \Om(M))^G
$$  
be the space of equivariant differential forms (resp. with
distributional coefficients.)  The derivation
$$ \d_G=1\otimes \d - v^a\otimes\iota_a$$
(using summation convention) of $C_G(M),\wh{C}_G(M)$ squares to $0$,
and defines the equivariant cohomologies $H_G(M),\
\wh{H}_G(M)$.  Both $H_G$ and $\wh{H}_G$ are functorial with respect
to maps between $G$-manifolds. In particular, the map $M\to \pt$ makes
$H_G(M)$ into a module over the ring $H_G(\pt)=(S\g^*_\C)^G$ of
invariant polynomials and $\wh{H}_G(M)$ into a module over the ring
$\wh{H}_G(\pt)=\ca{E}'(\g^*)^G$ of invariant compactly supported
distributions.

\subsubsection{$G$-valued equivariant cohomology}
In \cite{al:no} the following non-commutative Cartan complexes
$\ca{C}_G(M)$ and $\wh{\ca{C}}_G(M)$ were introduced.  Let
$\ca{E}'(G)$ be the convolution algebra of complex-valued
distributions on the group $G$. The subalgebra of distributions
supported at the group unit is isomorphic to the universal enveloping
algebra $U(\g)_\C$, with
$$ u_a=\frac{\d}{\d t} \ \delta_{\exp(te_a)} \Big|_{t=0}$$
representing the generator corresponding to the basis $e_a$.  Let
$$ \ca{C}_G(M)=(U(\g)_\C\otimes \Om(M))^G,\ \ \ \wh{\ca C}_G(M)
=(\ca{E}'(G)\otimes \Om(M))^G.$$
To define a differential on $\ca{C}_G(M), \wh{\ca C}_G(M)$ 
we assume that $\g$ comes equipped with an invariant inner product, 
and take the
basis $e_a$ to be orthonormal. One finds that
$$ {\d}_G=1\otimes\d-\f{u_a^L+u_a^R}{2}\otimes\iota_a +\f{1}{24}
f_{abc} 1\otimes \iota_a\iota_b\iota_c$$
(using summation convention) is a differential, i.e. $\d_G\circ\d_G=0$. 
Here $f_{abc}=[e_a,e_b]\cdot e_c$ are 
the structure constants for the given basis and
$u_a^L,u_a^R$ denote the operators of left-/right multiplication (that
is, convolution) by $u_a$.  We denote the corresponding cohomologies
by $\ca{H}_G(M),\wh{\H}_G(M)$. Notice that $\d_G$ is not a derivation
with respect to the obvious ring structure on $\wh{\ca
C}_G(M)$. However, one can modify the ring structure as follows: Let
$\Om(M)$ be equipped with the non-(super-)commutative ring structure
defined by
$$ \beta_1\odot\beta_2 =\diag_M^*\circ \exp(-\hh \iota_a^1\iota_a^2)
(\beta_1\times\beta_2)$$ where $\diag_M^*$ is pull-back by the
diagonal embedding, and $\iota_a^j$ are contractions acting on the
$j$th factor.  Together with the natural ring structure on
$U(\g),\,\ca{E}'(G)$, it defines a ring structure $\odot$ on ${\ca
C}_G(M),\,\wh{\ca C}_G(M)$ for which ${\d}_G$ is a derivation.  In
this way the cohomologies $\H_G(M)$ and $\wh{\H}_G(M)$ become
$\Z_2$-graded algebras. Both $\H_G$ and $\wh{\H}_G$ are functorial
with respect to maps of $G$-manifolds. In particular the map $M\to
\on{pt}$ makes $\H_G(M)$ into a module over the ring
$\H_G(\pt)=(U(\g)_\C) ^G$ of Casimir elements and $\wh{\H}_G(M)$ into
a module over the ring $\wh{\H}_G(\pt)=\ca{E}'(G)^G$ of invariant
distributions on the group.
	
\subsection{Weil models}
Below we will need not only the Cartan models of equivariant
cohomology but also the equivalent Weil models.
\subsubsection{$G$-differential algebras}
The concept of a $G$-differential algebra, introduced by Cartan
\cite{al:no}, generalizes the relations between contractions
$\iota_a$, Lie derivatives $L_a$, and exterior differential $\d$ on
the super-algebra $\Om^\star(M)$ of differential forms on a
$G$-manifold.

Define a graded vector-space $\hat{\g}= \bigoplus_{i\in\Z} \hat{\g}^i$
where $\hat{\g}^1=\R,\ \hat{\g}^{-1}=\hat{\g}^0=\g$ and
$\hat{\g}^j=\{0\}$ for $j\not= -1,0,1$.  Given $\xi\in\g$ let
$L_\xi\in\hat{\g}^0$ and $\iota_\xi\in\hat{\g}^{-1}$ be the
corresponding elements, and let $\d\in\hat{\g}^1$ be the generator
corresponding to $1\in\R$. Then $\hat{\g}$ is a graded Lie algebra
with brackets
$$ [\iota_\xi,\d]=L_\xi,\ \ [L_\xi,\iota_\eta]=\iota_{[\xi,\eta]},\ \
[L_\xi,L_\eta]=L_{[\xi,\eta]}$$
$$ [\iota_\xi,\iota_\eta]=0,\ \ [L_\xi,\d]=0,\ \ [\d,\d]=0$$
for all $\xi,\eta\in\g$.  Notice that $\hat{\g}^0\cong \g$ and
$\hat{\g}^{-1}\oplus \hat{\g}^0\cong \g\rtimes \g$ are Lie subalgebras
both of which are stable under $\d$.

A {\it $G$-differential space} is a super-vector space $B$, together
with a representation of $\hat{\g}$ on $B$. The subspace $B^G$
annihilated by $\hat{\g}^0$ is called the {\it invariant} subspace,
the subspace $B_{hor}$ annihilated by $\hat{\g}^{-1}$ is called the
{\it horizontal} subspace, and the subspace $B_{basic}$ annihilated by
$\hat{\g}^{-1}\oplus\hat{\g}^0$ is called the {\it basic}
subspace. Note that $B_{basic}$ is preserved under the differential
$\d$.

If $B$ carries in addition the structure of a super-algebra
\footnote{In this paper, (super)algebra always means 
(super)algebra with unit, and a homomorphism between (super)algebras
is required to take the unit element to the unit element.}, and if
$\hat{\g}$ is represented by derivations, then $B$ is called a
$G$-differential algebra.  If $B_1,B_2$ are $G$-differential spaces
then so is their ($\Z_2$-graded) tensor product $B_1\otimes B_2$.  A
{\it homomorphism of $G$-differential spaces} is a homomorphism of
super-spaces intertwining the $\hat{\g}$-representations.  Similar
definitions hold for $G$-differential algebras.

%
%


\subsubsection{The Weil algebras $W_G$ and $\wh{W}_G$.}
The super-commutative Weil algebras $W_G$ and $\wh{W}_G$ are defined
as tensor products
$$ W_G=S\g_\C^*\otimes \wedge\g_\C^*,\ \ \ \ \wh{W}_G=\ca{E}'(\g^*)\otimes 
\wedge\g_\C^*,$$
with $\Z_2$-grading inherited from the $\Z_2$-grading on 
the exterior algebra. They carry the structure of
$G$-differential algebras, as follows. Let $y^a\in\wedge\g^*$ 
be the generators corresponding to the basis $e^a$ of $\g^*$.
Let 
$$ L_a=L_a\otimes 1+1\otimes L_a$$
be the generators for the action on the Weil algebras $W_G,\wh{W}_G$
coming from the coadjoint action of $G$ on $\g^*$. 
The natural contraction operation 
on $\wedge\g^*$ extends to operators 
$\iota_a=1\otimes\iota_a$.
Let the Weil differential $\d^W$ be given by the formula 
$$
\d^W= y^a(L_a \otimes 1) + (v^a - \frac{1}{2} f^a_{bc} y^b y^c) \iota_a.
$$
%

Given a $G$-manifold $M$, the horizontal projection 
$P_{hor}:\,\wh{W}_G\to \ca{E}'(\g^*)$ induces  
isomorphisms between the basic subcomplexes and the Cartan algebras
$$ P_{hor}^W:=P_{hor}\otimes 1:\,
(\wh{W}_G\otimes\Om(M))_{basic}\to \wh{C}_G(M).$$
As shown in \cite{al:no}
this isomorphism takes $\d =\d^W\otimes 1+1\otimes \d$ to the 
Cartan differential $\d_G$, and it restricts to an 
ring isomorphism 
$P_{hor}^W:\,(W_G\otimes \Om(M))_{basic}\cong C_G(M)$. 
Therefore, 
$H_G(M)$ and $\wh{H}_G(M)$ are equivalently 
the cohomologies of $(W_G\otimes \Om(M))_{basic}$ and 
$(\wh{W}_G\otimes \Om(M))_{basic}$, respectively.

\subsubsection{Non-commutative Weil algebras}\label{subsec:noncom}
Let the Lie algebra $\g$ be equipped
with an invariant inner product, used to identify $\g\cong \g^*$, 
and suppose that the basis $e_a\in \g$
is orthonormal. 

Let $\Cl(\g)$ be the Clifford algebra of $\g$, 
defined as the quotient of the tensor algebra by the ideal generated 
by elements of the form  $2\xi\otimes\xi-\xi\cdot\xi$ with
$\xi\in\g$. The generators $x_a$ of $\Cl(\g)$ satisfy the 
super-bracket relations $[x_a, x_b]=\delta_{ab}$. 
The element $\gamma:=-\f{1}{6}f_{abc} x_a x_b x_c$ has the 
beautiful property \cite{ko:cl} that its square is a scalar: 
$\gamma^2=-\f{1}{48}f_{abc}f_{abc}$. 
Hence
$$ \d^{Cl}=\ad(-\f{1}{6}f_{abc} x_a x_b x_c):\,\Cl(\g)_\C\to\Cl(\g)_\C$$
is a differential. Together with contractions $\iota_a=\ad(x_a)$ and 
the Lie derivatives $L_a=\ad(-\hh f_{abc}x_b x_c)$  
for the adjoint action this makes $\Cl(\g)_\C$ into 
a $G$-differential algebra. 

In terms of the identification $\sig:\,\Cl(\g) \to\wedge\g$ 
by the symbol map, $\sigma (x_{a_1} \dots  x_{a_k}) =
y_{a_1} \dots  y_{a_k}$ for $a_1<\dots<a_k$,  
the differential $\d^{Cl}$ is given by the formula
\begin{equation}\label{CliffordDifferential} 
\d^{Cl}=
- \frac{1}{2} f_{abc}y_b y_c \iota_a - 
\frac{1}{24} f_{abc} \iota_a \iota_b \iota_c. 
\end{equation}

The Weil algebras $\W_G$ and $\wh{\W}_G$ are the non-commutative 
$G$-differential algebras 
$$
\W_G = U(\g)_\C  \otimes \Cl(\g)_\C \ , \ \wh{\W}_G= \ca{E}'(G) \otimes \Cl(\g)_\C
$$
with 
\beq \iota_a&=&\ad(x_a)\\
     L_a&=&\ad(u_a-\hh f_{abc}x_bx_c)\\  
     \d^\W&=& \ad(u_a x_a-\f{1}{6}f_{abc} x_a x_b x_c).
\eeq 
The operators $L_a$ generate the representation of $\g$ 
on the Weil algebra induced from the adjoint action on
$\g$.  

By \cite{al:no}, Proposition 3.7 we describe the Weil differential in
terms of the vector space isomorphism $\wh{\W}_G \cong
\ca{E}'(G)\otimes\wedge\g$ given by the symbol map:
\begin{equation} \label{sigmad}
\d^\W = y_a (L_a\otimes 1) + 
\left( \frac{u_a^L + u_a^R}{2} 
- \frac{1}{2} f_{abc}y_b y_c \right)
\iota_a - \frac{1}{24} f_{abc} \iota_a \iota_b \iota_c.
\end{equation}

Given a $G$-manifold $M$, the 
horizontal projection $P_{hor}:\,\wh{\W}_G\to \ca{E}'(G)$
induces isomorphisms of differential algebras, 
$$ P_{hor}^\W=P_{hor}\otimes 1:\,
(\wh{\W}_G\otimes\Om(M))_{basic}\cong \wh{\ca{C}}_G(M),
 \ \ \ \ ({\W}_G\otimes\Om(M))_{basic}\cong {\ca{C}}_G(M)$$
so that the equivariant cohomologies $\H_G(M)$ and 
$\wh{\H}_G(M)$ are equivalently the cohomologies 
of $({\W}_G\otimes\Om(M))_{basic}$ and $(\wh{\W}_G\otimes\Om(M))_{basic}$. 

\subsubsection{The quantization map}
Let $J^\hh\in C^\infty(\g)^G$ be the unique smooth square root of the
determinant of the Jacobian of the exponential map such that
$J^\hh(0)=1$. The {\em Duflo map} \cite{du:op} is the linear map
\begin{equation}\label{DufloMap}
 \on{Duf}:\,\ca{E}'(\g)\to \ca{E}'(G),\ \ \ v\mapsto \exp_*(J^\hh v).
\end{equation}
It takes distributions supported at the origin to distributions
supported at the group unit, hence defines a vector space isomorphism
$\on{Duf}:\,S(\g)\to U(\g)$. The important property of the Duflo map
is that it becomes a {\em ring homomorphism} if restricted to
invariant elements.

In \cite{al:no} we define an explicit homomorphism of $G$-differential
spaces, called the quantization map,
$$ \ca{Q}:\,\wh{W}_G\to \wh{\W}_G $$
which on $W_G$ restricts to an isomorphism, $W_G\to \W_G$.  On
$1\otimes \wedge\g^*_\C$ the map $\ca{Q}$ coincides with the symbol
map, and on $\ca{E}'(\g) \otimes 1$ it is equal to the Duflo
map. Moreover, we show that the induced homomorphisms in cohomology
preserve the ring structures, i.e. the map $\wh{H}_G(M)\to
\wh{\H}_G(M)$ is an algebra homomorphism and the map $H_G(M)\to
\H_G(M)$ is an algebra isomorphism.

\section{Abelianization}
\label{Sec:Abel}
The equivariant cohomology groups $H_G(M)$ and $\wh{H}_G(M)$ are
functorial not only with respect to maps of $G$-manifolds but also
with respect to group homomorphisms. In particular, the inclusion of
the maximal torus $T$ induces maps $R_G^T:\,H_G(M)\to H_T(M)$ and
$R_G^T:\, \wh{H}_G(M)\to \wh{H}_T(M)$. Both of these maps are
injective, and their image is the subalgebra of Weyl invariants.

This functoriality property does not carry over to group-valued
equivariant cohomology $\wh{\H}_G(M)$, because a group homomorphism
$\phi:\,K\to G$ does not give rise to $K$-equivariant maps ${\W}_G\to
{\W}_K$ or $\wh{\W}_G\to \wh{\W}_K$ in any natural way.

It is the goal of this section to obtain a suitable replacement for
the restriction map $R_G^T$.  The main steps are as follows.
\begin{enumerate}
\item 
First, we define a new Weil algebra $\wh{\W}_T^-=\ca{E}'(T)^-\otimes
\Cl(\t)_\C$,
involving the negative Hardy space $\ca{E}'(T)^-$. 
\item 
We construct a map $r_G^T:\,\ca{E}'(G)\to \ca{E}'(T)^-$ with nice 
properties. 
\item
We define a map $p_G^T\circ \ca{T}:\,\Cl(\g)_\C\to\Cl(\t)_\C$, 
where $p_G^T$ is orthogonal projection and 
the ``twist''  operator $\ca{\T}$ 
involves the classical r-matrix. We show that this 
map intertwines the Clifford differentials ``up to a $\rho$-shift''. 
\item
The restriction map $\ca{R}_G^T:\,\wh{\W}_G\to\wh{\W}_T^-$ will 
be defined as a product of the two maps, 
$r_G^T$ and $p_G^T\circ\ca{T}$. 
The main result of this Section is Theorem \ref{restriction}, 
stating that $\ca{R}_G^T$ is a homomorphism of 
$T$-differential spaces. It therefore defines a map in equivariant 
cohomology.
\end{enumerate}

Throughout this section we will assume that $G$ is a product 
of a simply connected group and a torus (although only 
the simply connected part will play a role). 
A given invariant inner product on $\g$ will be used to identify 
$\g\cong\g^*$.

We let $T$ be a maximal torus of $G$, with Lie algebra $\t$, 
and denote by $\t_+$ some choice of a fundamental 
Weyl chamber. Let $\Lambda\subset\t$ be the integral lattice, 
its dual $\Lambda^*\subset\t $ the weight lattice, and 
$\Lambda^*_+=\Lambda^*\cap\t_+$ the set of dominant weights.

Given $t\in T$ and $\lambda\in\Lambda^*$ we denote 
by $t^\lambda\in S^1\subset\C$ 
the image of $t$ under the $T$-character 
$T\to S^1$ defined by $\lambda$, that is 
$(\exp\xi)^\lambda=e^{\tpi\l\lambda,\xi\r}$. 
Abusing notation we sometimes denote the map itself by 
$t^\lambda$. 

Let $\mf{R}\subset\Lambda^*$ be the set of roots of $G$ and 
$\mf{R}_+$ the subset of positive roots.  
Our assumption that $G$ be a product of a connected, simply 
connected Lie group and a torus, implies that the 
half-sum of positive roots 
$$ \rho=\hh \sum_{\alpha\in\mf{R}_+}\alpha$$
is a weight $\rho\in\Lambda^*_+$. It is alternatively characterized
as the sum of fundamental weights for the semi-simple part 
of $G$.

For every positive root $\alpha\in\mathfrak{R}_+$ 
let $e_\alpha\in\g_\C$ be a root 
vector, with property $[\xi,e_\alpha]=\tpi\l\alpha,\xi\r 
e_\alpha$ for $\xi\in\t$, and set $e_{-\alpha}=\ol{e_\alpha}$.
Normalize the length of $e_\alpha$ so that $e_{\alpha}\cdot
e_{-\alpha}=1$. Then 
\begin{equation}\label{Normalization}
 [e_\alpha,e_{-\alpha}]=\tpi \alpha.
\end{equation}

Every weight $\lambda\in\Lambda^*_+$ 
parametrizes a unique irreducible $G$-representation $V_\lambda$, 
containing a highest weight vector $v_\lambda\not=0$, 
with 
\begin{equation}\label{PropertiesHighestWeight}
e_\alpha\cdot v_\lambda=0,\ \ \ 
\end{equation}
for every positive root $\alpha\in \mathfrak{R}_+$.
and 
\begin{equation}\label{PropertiesHighestWeight2}
t\cdot v_\lambda=t^\lambda v_\lambda
\end{equation}
for every $t\in T$. We denote by $\chi_\lambda$ 
the character of the representation $V_\lambda$.

\subsection{The Hardy-Weil algebra $\wh{\W}_T^-$}
Recall (cf. \cite{fr:in}, Chapter 8.5) 
that for any distribution  $u\in\ca{E}'(T)$, the 
Fourier coefficients $\l u,t^{\lambda}\r$ 
are polynomially bounded: 
$$|\l u, t^\lambda\r\,|\le C(1+|\lambda|)^N$$
for some $C>0, N\in \Z$. 
Conversely every sequence of complex numbers 
satisfying such an estimate defines a distribution on $T$.
We can therefore identify $\ca{E}'(T)$ with the space of 
polynomially bounded sequences $\Lambda^*\to\C$. 
Similarly, if $\phi\in C^\infty(T)$ has Fourier decomposition
$\phi=\sum_{\nu\in\Lambda^*}a_\nu t^\nu$ then the coefficients 
$a_\nu$ are rapidly decreasing, that is for all $k>0$ there exists 
$C_k$ with 
$$ |a_\nu|\le C_k (1+|\nu|)^{-k}.$$
Conversely, every rapidly decreasing sequence $a_\nu$ defines 
a smooth function on the torus.

Let the negative Hardy space
$\ca{E}'(T)^-\subset\ca{E}'(T)$ consist  
of all distributions on $T$ for which only strictly 
negative Fourier coefficients are non-vanishing,
$$
\ca{E}'(T)^-:= \{u \in \ca{E}'(T)|\l u,t^\nu \r =0 \
{\rm if}\  \nu \not\in (\rho+\Lambda^*_+) \},
$$
and let $\pi^-$ be the Szeg\"o projector onto this space, 
$$ \l \pi^-(u),t^{\lambda+\rho}\r=\left\{ 
\begin{array}{cl} 
\l u,t^\nu\r &\mbox{ if } \lambda \in \Lambda^*_+,\\
0&\mbox{ otherwise .}
\end{array}\right.
$$
Convolution of distributions makes the Hardy space into an 
algebra with unit element $ \delta^-_e:=\pi^-(\delta_e)$. 
The Szeg\"o projector is an algebra homomorphism. (The inclusion of 
$\ca{E}'(T)^-$ into $\ca{E}'(T)$ is not an algebra homomorphism 
since the unit elements are different.) 

The Hardy-Weil algebra is defined as the tensor product 
$$ \wh{\W}_T^-:=\ca{E}'(T)^-\otimes \Cl(\t)_\C.$$
It is a $T$-differential algebra, with differential 
$\d^-:=\sum_{k=1}^r u_k^-\iota_k$ where $u_k^-=\pi^-(u_k)$. 
The projection $\pi^-:=\pi^-\otimes 1:\,
\wh{\W}_T\to \wh{\W}_T^-$ is a homomorphism of $T$-differential 
algebras 
and the inclusion $\wh{\W}_T^-\hra \wh{W}_T$ is 
a homomorphism of $T$-differential spaces (not algebras). 

Given a $T$-manifold $M$, let $\wh{\H}_T^-(M)$ be the cohomology of the 
basic subcomplex $(\wh{\W}_T^-\otimes \Om(M))_{basic}$. 
Equivalently, it is the cohomology for the Cartan model, 
$$\wh{\ca{C}}_T^-(M)=\ca{E}'(T)^-\otimes \Om(M)^T,\ \ 
\d_T^-:=1\otimes \d-\sum_{k=1}^r u_k^-\otimes \iota_k.$$ 
The projection $\wh{\W}_T\to \wh{\W}_T^-$
gives rise to a ring homomorphism 
$ \wh{\H}_T(M)\to \wh{\H}^-_T(M)$, 
and the inclusion $\wh{\W}_T^-\hra \wh{W}_T$ 
defines a linear map $  \wh{\H}^-_T(M)\to \wh{\H}_T(M)$. 

\subsection{The restriction map $r_G^T:\,\ca{E}'(G)\to \ca{E}'(T)^-$}
Below we will often find it convenient to work with a 
Cartan-Weyl basis in $\g_\C$, given by the collection of 
$e_\alpha,e_{-\alpha}$ supplemented by a basis
$e_k$ ($k=1,\ldots,r$) for $\t$.

For all $\lambda\in\Lambda^*_+$ let $v_\lambda$ be a highest 
weight vector of unit length, and define a 
function $\S_\lambda\in C^\infty(G)$ as the matrix element 
\begin{equation}\label{DeltaExplicit}
\S_\lambda(g)=( v_\lambda,g\cdot
v_\lambda),
\end{equation}
using the inner product $(\cdot,\cdot)$ on $V_\lambda$. 
It has the properties
\begin{equation}
\S_\lambda(e)=1,
\end{equation}
\begin{equation}
\S_\lambda(tg)=\S_\lambda(gt)=t^\lambda\S_\lambda(g)
,\ \ (t\in T,\,g\in G),
\end{equation}
\begin{equation} \label{evanish}
e_{\alpha}^L \S_\lambda=e_{-\alpha}^R \S_\lambda =0.
\end{equation}
for $\alpha\in \mf{R}_+$, where for any $\xi\in\g_\C$ we denote by 
$\xi^L,\xi^R$ the left-and right invariant vector fields
equal to $\xi$ at the group unit. 
In fact, because of the decomposition 
$$L^2(G)=\bigoplus_{\lambda\in\Lambda^*_+} V_\lambda\otimes 
V_{\lambda}^*=\bigoplus_{\lambda\in\Lambda^*_+} 
\on{End}(V_\lambda)
$$ 
as a $G\times G$ representation,  
the function $\S_\lambda$ is uniquely 
characterized by these three properties.
The functions $\S_\lambda$ have the property 
\begin{equation}\label{Multaplicativity}
\S_{\lambda + \nu}= \S_{\lambda} \S_{\nu}
\end{equation}
for all $\lambda,\nu\in\Lambda^*_+$. 
\begin{example}
Let $G = SU(2)$. All dominant weights are multiples $\lambda=k\rho$ 
with $k=0,1,\ldots$, and 
$$
\S_{k\rho} \left( 
\begin{array}{cc}
a & b \\
c & d 
\end{array}
\right)  =a^k.
$$ 
\end{example}

\begin{proposition}
The map $r_G^T:\,\ca{E}'(G)\to \ca{E}'(T)^-$, 
$$ 
\l r_G^T(u),t^{\lambda+\rho}\r=\l u,\S_\lambda\r 
$$ 
is well-defined and continuous. 
\end{proposition}

\begin{proof}
To show that $\l u,\S_\lambda\r$ are Fourier coefficients for 
a well-defined distribution $r_G^T(u)\in\ca{E}'(T)^-$, we have 
to show that they are polynomially bounded as a function of 
$\lambda$. Clearly, this only involves the semi-simple 
part of $G$ (indeed, if $G$ is abelian, there is nothing to prove since 
$\ca{E}'(T)^- = \ca{E}'(G)$ in that case and the map $r_G^T$ is
the identity map). 
 
Assume therefore that $G$ is simply connected. 
Then $\Lambda^*_+$ is spanned by fundamental weights 
$w_1,\ldots,w_r$ of $G$. Let 
$$\S:= (\S_{w_1}, \dots , \S_{w_r}): G \to \C^r.$$
By \eqref{Multaplicativity}, $\S_\lambda= \S ^* z^\lambda$ using
multi-index notation $z^{\lambda}=\prod_j z_j^{\lambda_j}$, hence $\l
u,\S_\lambda\r=\l\S_* u,z^\lambda\r$.  The image of $\S$ is contained
in the unit polydisk $D=\{ (z_1, \dots, z_r)\in\C^r| \ |z_k| \le 1
\}\subset\C^r$.

For any distribution $v\in\ca{E}'(\C^r)$ supported on the polydisk
$D$, there exist constants $C>0$ and $k\in\N$ such that the value of
$v$ on any test function $\phi$ is bounded by derivatives of $\phi$ up
to order $k$ (cf. H\"ormander \cite{ho:an}, Theorem 2.3.10)
$$ |(v,\phi)|\le C \sum_{|J|\le k}{\sup}_{x\in D}\Big|
\f{\p^{|J|}\phi}{\p x^J}\Big|.$$
Applying this to $v=\S_*u$ and $\phi=z^\lambda$ immediately gives 
the required estimates.

To show that $r_G^T$ is continuous, consider
the semi-norms $\l r_G^T(u),\phi\r$ for $\phi\in C^\infty(T)$. 
Since the coefficients $a_\nu$ in the Fourier decomposition 
$$\phi=\sum_{\nu\in\Lambda^*} a_\nu t^\nu\in C^\infty(T)$$
are rapidly decreasing, and since $|\S_\lambda|\le 1$, the function 
$$ \psi:=\sum_{\lambda\in\Lambda^*_+} a_{\lambda+\rho}\S_\lambda $$
on $G$ is well-defined and continuous. 
Let $B$ be the Laplace operator on $G$. On $V_\lambda\otimes
V_\lambda^*\subset L^2(G)$ 
it acts as a scalar $4\pi^2(|\lambda+\rho|^2-|\rho|^2)$.
Hence 
$$ B^k\psi=\sum_{\lambda\in\Lambda^*_+} 
(4\pi^2(|\lambda+\rho|^2-|\rho|^2))^k 
a_{\lambda+\rho}\S_\lambda 
$$
is still continuous. By the Sobolev Lemma, 
this shows that $\psi$ is smooth. Since 
$ \l r_G^T(u),\phi\r=\l u,\psi\r$
it follows that $r_G^T$ is continuous. 
\end{proof}

\begin{example}\label{deltaminus}
Given $h\in T$, let $\delta_h^-\in\ca{E}'(T)^-$ be the distribution 
$$ \delta_h^-=\pi^-(\delta_h^T),$$
where $\delta_h^T$ is the delta-distribution on $T$. 
Note that  
$$ \delta_{h_1}^-\delta_{h_2}^-=\delta_{h_1\,h_2}^-$$
for $h_1,h_2\in T$. The distribution $\delta_h^-$ on $T$ 
is related to the delta-distribution 
$\delta_h$ on $G$ by 
$$ h^\rho r_G^T(\delta_h)=\delta_h^-,$$
by the calculation 
$\l r_G^T(\delta_h),t^{\lambda+\rho}\r=\Delta_\lambda(h)=h^\lambda$.
\end{example}

\begin{proposition}[Properties of the map $r_G^T$]
Let $u_{-\alpha}\,,u_k\,,u_{\alpha}$ be the generators of 
$U(\g)_\C$ in the Cartan-Weyl basis. Then 
\begin{eqnarray}
\label{eq:8}
r_G^T\circ (u_k)^R=r_G^T\circ (u_k)^L&=&(u_k^-
-\tpi\rho_k)\circ r_G^T\\
\label{eq:9}
r_G^T\circ (u_{-\alpha})^L=r_G^T\circ (u_\alpha)^R&=&0
\end{eqnarray}
for $k=1,\ldots,r$ and all $\alpha\in\mathfrak{R}_+$.
Here the superscripts $L/R$ denote multiplication 
from the left/right. 
\end{proposition}

\begin{proof}
%
%
Let $u\in\ca{E}'(G)$ and $\lambda\in\Lambda^*_+$ be given. Then 
$$ \l u_k^-\circ r_G^T(u),t^{\lambda+\rho}\r =
\tpi(\lambda_k+\rho_k) \l r_G^T(u),t^{\lambda+\rho}\r$$
while on the other hand
$$ \l r_G^T(u_k u),\,t^{\lambda+\rho}\r
=\l (u_k u), \S_\lambda\r
=\l u, e_k^R \S_\lambda\r
=\tpi \lambda_k \l u,\S_\lambda\r
=\tpi \lambda_k \l r_G^T(u),\,t^{\lambda+\rho}\r.
$$
Hence, $r_G^T\circ (u_k)^L=(u_k^- -\tpi\rho_k)\circ r_G^T$. 
Similarly, $r_G^T\circ (u_k)^R =(u_k^- -\tpi\rho_k)\circ r_G^T$, 
proving  \eqref{eq:8}. 
The equations \eqref{eq:9} are obtained from 
\eqref{evanish} as follows:
$$ \l r_G^T(u_{-\alpha}u ),t^{\lambda+\rho}\r
=\l u_{-\alpha}u ,\S_\lambda\r=
-\l u ,e_{-\alpha}^R\S_\lambda\r=0,$$
$$ \l r_G^T(u  u_{\alpha}),t^{\lambda+\rho}\r
=\l u  u_{\alpha},\S_\lambda\r=
\l u ,e_{\alpha}^L\S_\lambda\r=0.$$
\end{proof}

The function $\S_\lambda$ is $\Ad(T)$-invariant but not 
$\Ad(G)$-invariant. Let $P^G:\,C^\infty(G)\to C^\infty(G)^G$ be  
projection to the $\Ad(G)$-invariant part (given by 
averaging over the group). We claim that
\begin{equation}\label{inv}
P^G(\S_\lambda)=(\dim V_\lambda)^{-1}\chi_\lambda. 
\end{equation}
To see this let $\tau_\lambda:\, G\to \on{Aut}(V_\lambda)$ be 
the representation labeled by $\lambda$. 
For any weight $\nu\in \Lambda^*_+$, the 
operator $\int_G \tau_\lambda(g) \ol{\chi_\nu(g)} \d g$ 
commutes with all $\tau_\lambda(h)$ and is therefore a multiple 
of the identity:
$$ \int_G \tau_\lambda(g) \ol{\chi_\nu(g)} \d g 
=\delta_{\nu,\lambda}\f{\Vol G}{\dim V_\lambda} 
\on{Id}_{V_\lambda},$$
where the constant is verified by taking the trace of both sides.
As a consequence, 
$$ 
\int_G \Delta_\lambda(g) \ol{\chi_\nu(g)} \d g =
\int_G (v_\lambda ,g\cdot v_\lambda) \ol{\chi_\nu(g)} \d g  
=\delta_{\nu,\lambda}\f{\Vol G}{\dim V_\lambda}
$$
which implies \eqref{inv}.

\subsection{Restriction for $\Cl(\g)_\C$}
As our next ingredient in the definition of the map
$\ca{R}_G^T:\,\wh{\W}_G\to\wh{\W}_T^-$ we need a suitable projection
from $\Cl(\g)_\C\to \Cl(\t)_\C$.  Let $\Cl(\g)_\C=\Cl(\t)_\C\otimes
\Cl(\t^\perp)_\C$ be the decomposition of the Clifford algebra, and
$p_G^T:\,\Cl(\g)_\C\to
\Cl(\t)_\C$ the orthogonal projection defined by it.  
We will pre-compose this projection with a ``twist''
$$\ca{T}=\exp(\hh r_{ab}\iota_a\iota_b)$$ 
where $r = r_{ab}e_a e_b
\in \wedge^2 \g_\C$ is the classical r-matrix,
$$ r=\sum_\alpha e_\alpha\wedge e_{-\alpha}
=
\hh \sum_\alpha (e_\alpha\otimes e_{-\alpha}
-e_{-\alpha} \otimes  e_\alpha),$$
Recall the following properties of the r-matrix:
\begin{proposition}[Properties of the classical r-matrix] 
\label{Rproperties}
The $r$-matrix $r\in \wedge^2 \g_\C$ satisfies the equations 
\begin{enumerate}
\item (Classical Yang-Baxter Equation) 
\begin{equation}\label{YBE}
 \on{Cycl}_{abc}\big( r_{as}f_{sbt}r_{tc}\big)
=\f{1}{4} f_{abc} 
\end{equation}
where $\on{Cycl}_{abc}$ is the sum over cyclic permutations.
\item
\begin{equation}\label{rrho}
f_{abc}r_{bc} =4\pi i \rho_a.
\end{equation}
\item
For any $\mu\in\t$, 
\begin{equation}\label{invt}
\mu_a (f_{abs}r_{sc}-f_{asc}r_{bs})=0.
\end{equation}
\end{enumerate} 
\end{proposition}
We include the proof of Proposition \ref{Rproperties}, 
since many different normalizations of the 
$r$-matrix appear in the literature.
\begin{proof}
Using the Cartan-Weyl basis for $\g_\C$, the only 
non-vanishing entries in $r_{ab}$ are $r_{-\alpha, \alpha}=
-r_{\alpha,-\alpha}=1/2$. Up to permutation of indices, 
the only non-vanishing structure
constants are of the form
$f_{-\alpha, \alpha, i}\, , f_{-\alpha - \beta, \alpha, \beta}$. 
One verifies:
\beq r_{-\alpha, \alpha}\, f_{-\alpha, i, \alpha}\, r_{-\alpha,
\alpha}& =& \frac{1}{4} f_{-\alpha, i, \alpha} , \\ r_{-\alpha,
\alpha}\, f_{-\alpha, -\beta, \alpha+\beta} \,
r_{-\alpha-\beta,\alpha+\beta} + r_{-\beta, \beta}\, f_{-\beta,
\alpha+\beta, -\alpha}\, r_{\alpha, -\alpha} \\+
r_{\alpha+\beta,-\alpha-\beta\,} f_{ \alpha+\beta, -\alpha, -\beta}\,
r_{\beta, -\beta} &=& \frac{1}{4} f_{-\alpha, -\beta, \alpha + \beta}.
\eeq Equation \eqref{rrho} follows from the calculation, using
\eqref{Normalization},
%
$$ r_{ab} f_{abc} e_c =\sum_{\alpha\in\mathfrak{R}_+}
[e_\alpha,e_{-\alpha}]=\tpi \sum_{\alpha\in\mathfrak{R}_+}
\alpha=4\pi i \rho.
$$
The last Equation \eqref{invt} is just the infinitesimal version of 
the $T$-invariance of $r$.
\end{proof}
The main property of the composition $p_G^T\circ \ca{T}$ 
is that it intertwines the Clifford differential on 
$\Cl(\g)_\C$ with a very simple differential on $\Cl(\t)_\C$:
\begin{proposition}\label{Clif}
The composition $p_G^T\circ \ca{T}$ intertwines the Clifford 
differential
$\d^{Cl}$ on $\Cl(\g)_\C$ and the differential $\tpi\sum_{k=1}^r
\rho_k\iota_k$ on $\Cl(\t)_\C$:
$$ (p_G^T\circ \ca{T})\circ \d^{Cl}=
(\tpi \sum_{k=1}^r \rho_k \iota_k)\circ (p_G^T\circ \ca{T}).$$
\end{proposition}

\begin{proof}
We use the symbol map $\sig:\,\Cl(\g)\to\wedge\g$ 
to identify $\d^{Cl}$ with the 
differential \eqref{CliffordDifferential} on $\wedge\g_\C$. 
Since $\Ad(\T) y_a =y_a - r_{ab} \iota_b$, we find
$$
\Ad(\T) \d^{Cl}=
\left(- \frac{1}{2} f_{abc} (y_b - r_{bl}\iota_l)
(y_c - r_{cm} \iota_m) \right) \iota_a
  - \frac{1}{24} f_{abc} \iota_a \iota_b \iota_c.
$$
In this expression, 
the terms cubic in contractions $\iota_a$ cancel thanks to the Yang-Baxter
equation \eqref{YBE} 
for $r$. The remaining terms 
combine, using \eqref{rrho}, to 
\beq 
\Ad(\T)\d^{Cl}&=&-\hh f_{abc} y_b y_c \iota_a+
\hh f_{abc}(r_{cm}y_b \iota_m\iota_a+r_{bl}\iota_l y_c \iota_a
)\\
&=& -\hh f_{abc} y_b y_c \iota_a+
f_{abc}(r_{cm} y_b\iota_m\iota_a+ \hh r_{bc}\iota_a)  \\
&=& -\hh f_{abc} y_b y_c \iota_a 
+f_{abd}r_{dc} y_a\iota_b\iota_c+\tpi \rho_a\iota_a.
\eeq
Composing this operator with $p_G^T$ 
kills the first term, and also the second since 
$$ p_G^T\circ (y_a f_{abd} r_{dc}\iota_b\iota_c)=
\sum_{k=1}^r y_k\circ 
p_G^T\circ f_{kbd}r_{dc}\iota_b\iota_c=0
$$
where we have used Equation \eqref{invt}. We conclude 
$$ p_G^T \circ \Ad(\T)\d^{Cl}=\tpi p_G^T\circ
\sum_{k=1}^r \rho_k\iota_k  
=\tpi \sum_{k=1}^r \rho_k\iota_k 
\circ p_G^T.$$
\end{proof}

\subsection{Restriction for $\wh{\W}_G$}
Taking the product of the restriction maps 
$r_G^T:\,\ca{E}'(G)\to\ca{E}'(T)^-$ and 
$(p_G^T \circ \T):\,\Cl(\g)_\C\to\Cl(\t)_\C$
we obtain a map between Weil algebras, 
$$\ca{R}_G^T:= r_G^T\otimes (p_G^T \circ \T)
: \wh{\W}_G \to \wh{\W}_{T}^-.$$
The main result of this Section reads:
\begin{theorem}  \label{restriction}
The map $\ca{R}_G^T: \wh{\W}_G \to \wh{\W}_{T}^-$ is 
a homomorphism of $T$-differential spaces. That is, it
is a chain map which moreover 
intertwines the contractions and 
Lie derivatives for the $T$-action.
\end{theorem}

\begin{proof}
It is clear that $\ca{R}_G^T$ commutes with Lie derivatives 
$L_\xi$ and contractions $\iota_\xi$ for all $\xi\in\t$. 
To show that $\ca{R}_G^T$ is a chain map 
we show that $(r_G^T\otimes p_G^T)\circ \Ad(\T)\d^\W
=\d^\W \circ (r_G^T\otimes p_G^T)$. By definition, the 
Weil differential splits into two pieces
$\d^\W=\d'+ \d^{Cl}$, with $\d'=\ad(u_a x_a)$. 
Again we identify $\Cl(\g)\cong \wedge\g$ by the symbol map. 
Since $x_a^L=y_a-\hh\iota_a$ and 
$x_a^R=y_a+\hh\iota_a$, we have: 
$$\d'= y_a (u_a^L-u_a^R)+\hh (u_a^L+u_a^R)\iota_a.$$ 
Conjugating by $\T$, and using  $ \Ad(\T)y_a=y_a-r_{ab}\iota_b$, 
$$ \Ad(\T)\d'=(y_a-r_{ab}\iota_b) (u_a^L-u_a^R)+\hh (u_a^L+u_a^R)\iota_a.$$ 
In the Cartan-Weyl basis,
$$
r_{ab}\iota_b(u_a^L-u_a^R) =\hh \sum_{\alpha\in\mathfrak{R}_+}
\big((u_\alpha^L-u_\alpha^R)\iota_{-\alpha}-
(u_{-\alpha}^L-u_{-\alpha}^R)\iota_{\alpha}\big)
$$
while
$$\hh(u_a^L+u_a^R)\iota_a=\hh \sum_{k=1}^r (u_k^L+u_k^R)\iota_k
+\hh \sum_{\alpha\in\mf{R}_+}\big((u_\alpha^L+u_\alpha^R)\iota_{-\alpha}
+(u_{-\alpha}^L+u_{-\alpha}^R)\iota_{\alpha}\big).
$$
Hence, 
$$
 \Ad(\T)\d'=L_a\otimes y_a
+\sum_{\alpha\in\mathfrak{R}_+}(u^L_{-\alpha}\iota_\alpha
+u_\alpha^R\iota_{-\alpha})
+\sum_{k=1}^r \left(\f{u_k^L+u_k^R}{2}\right)\iota_k.
$$
Next, apply $r_G^T\otimes p_G^T$ to this expression. 
This kills the first term since 
$$ (r_G^T\otimes p_G^T)(L_a\otimes y_a)=
\sum_{k=1}^r (r_G^T\circ L_k)\otimes y_k=0,$$ 
(by $\Ad(T)$-invariance of $r_G^T$). 
Since $r_G^T\circ u_{-\alpha}^L=r_G^T \circ u_{\alpha}^R=0$,
the second term is killed as well. Only the last term survives 
and using 
\eqref{eq:8} we find
\beq (r_G^T\otimes p_G^T)\circ \Ad(\T)\d'&=&(r_G^T\otimes p_G^T)
\sum_{k=1}^r \left(\f{u_k^L+u_k^R}{2}\right)\iota_k
\\&=&
\sum_{k=1}^r (u_k^- -\tpi \rho_k)\iota_k \circ(r_G^T\otimes p_G^T).
\eeq
Together with Proposition \ref{Clif} this shows
$$ (r_G^T\otimes p_G^T)\circ \Ad(\T)(\d'+d^{Cl})=
\sum_{k=1}^r u_k^- \iota_k \circ(r_G^T\otimes p_G^T)
$$
which completes the proof.
\end{proof}

\subsection{Restriction $\ca{R}_G^T$ in cohomology}
For every $G$-manifold $M$, the homomorphism of $T$-differential 
algebras $\ca{R}_G^T:\,\wh{\W}_G\to \wh{\W}_T^-$ 
studied in the previous section induces a chain map 
between the basic subcomplexes
$$ \ca{R}_G^T=\ca{R}_G^T\otimes 1:\,(\wh{\W}_G\otimes\Om(M))_{basic}\to
(\wh{\W}_T^-\otimes\Om(M))_{basic},
$$
hence in cohomology, 
$$ \ca{R}_G^T: \wh{\H}_G(M) \to \wh{\H}_{T}^-(M).$$
Using the horizontal projection $P_{hor}^\W$ 
to identify Weil and Cartan models, the restriction map 
gives rise to a map $\wh{\ca{C}}_G(M)\to \wh{\ca{C}}_T^-(M)$. 
This map is described in the following Proposition: 
\begin{proposition}
The map between Cartan models induced by the restriction 
$\ca{R}_G^T$ is the map
$$ r_G^T\otimes \exp(\hh\iota(r_M)):\,\wh{\ca{C}}_G(M)\to
\wh{\ca{C}}_T^-(M).
$$
Here 
$r_M=r_{ab}(e_a)_M \wedge (e_b)_M$ 
is the bi-vector field on $M$ corresponding to the classical 
$r$-matrix.  
\end{proposition}

\begin{proof}
We have to compute the composition
$$P_{hor}^\W\circ (r_G^Tp_G^T\otimes 1)\circ
\exp(\hh r_{ab}\iota_a \iota_b\otimes 1):\,(\wh{\W}_G\otimes \Om(M))_{basic} 
\to \wh{\ca{C}}_G(M).
$$
Since $\iota_a\otimes 1=-1\otimes \iota_a$ 
on basic elements, the operator $\exp(\hh
r_{ab}\iota_a\iota_b\otimes 1)$ can be replaced with 
$\exp(\hh r_{ab}
(1 \otimes 
\iota_a\iota_b))=\exp(\hh\iota(r_M))$.
The latter commutes with $P_{hor}^\W\circ (r_G^Tp_G^T\otimes 1)$, 
and since $P_{hor}^\W\circ r_G^Tp_G^T=
r_G^T\circ P_{hor}^\W$ it follows that the above composition equals
$$  \big(r_G^T\otimes \exp(\hh \iota(r_M))\big)
\circ P_{hor}^\W $$ 
as claimed.
\end{proof}

\section{Properties of the restriction map}
\label{Sec:PropRes}
Having defined the restriction map in cohomology, some
natural questions arise: Does restriction commute with 
quantization? Is restriction a ring map? 
In this Section we will give an affirmative answer 
to the first question and a partial answer to the second. 
The strategy for approaching these problems is to describe 
the relevant maps in terms of their integral kernels. 
Section \ref{Sec:Duality1} and \ref{Sec:Duality2} review 
material from \cite{al:no}. 

\subsection{Duality for $\wh{W}_G$}\label{Sec:Duality1}
The Weil differential simplifies under the isomorphism of $\wh{W}_G$ 
given by multiplication with the function 
$\tau_0\in C^\infty(\g^*)\otimes \wedge \g^*$, 
$$
\tau_0(\mu) = \exp\left( - \frac{1}{2} f_{ab}^c y^a y^b \mu_c \right). 
$$
Indeed, 
\begin{equation}\label{ConjWeil}
 \Ad(\tau_0^{-1})\d^W=v^a\iota_a,\ \ \ 
\Ad(\tau_0^{-1})\iota_a=\iota_a-f_{ab}^c\mu_c y^b.
\end{equation}
These equations exhibit the Weil differential and the contraction as
dual to certain operators on the space $\Om(\g^*)$ of differential
forms on $\g^*$.  Indeed let $\bfl\cdot,\cdot\bfr$ denote the natural
pairing between $\wh{W}_G=\ca{E}'(\g^*)\otimes
\wedge\g^*$ and the space
$\Om(\g^*)=C^\infty(\g^*)\otimes\wedge\g$, and define a new pairing by 
$\l A,\beta\r:=\bfl\tau_0^{-1}A,\beta\bfr$. We denote by 
$\d^{Rh}$ the de Rham differential on $\Om(\g^*)$
and by $\iota_a^{Rh}$ the contractions with respect to the 
conjugation action.

Let $(\cdot)^t$ 
denote the transpose with respect to the pairing $\l\cdot,\cdot\r$. 
Then Equations \eqref{ConjWeil} become 
\begin{equation}\label{CW2}
\d^W=-(\d^{Rh})^t,\ \ \ \iota_a=-(\iota_a^{Rh}+\d\mu_a)^t
\end{equation}
where the differential $\d\mu_a$ acts by multiplication on $\Om(\g^*)$.
 
The integral kernel of the identity map $\wh{W}_G\to \wh{W}_G$ with
respect to the pairing $\l\cdot,\cdot\r$ is the element
\begin{equation}\label{Def:Lambda0}
\Lambda_0 := e^{-y^a d \mu_a} \tau_0(\mu) \delta_{\mu}
\in(\wh{W}_G \otimes \Omega(\g^*))^G.
\end{equation}
That is, $ \l A, \Lambda_0 \r = A $ for any $A \in \wh{W}_G$.  It
follows directly from \eqref{CW2} that
\begin{equation} \label{Lambda0}
(\d^W\otimes 1+1\otimes \d^{Rh}) \Lambda_0 =0, \ 
(\iota_a\otimes 1+1\otimes \iota_a^{Rh}) 
\Lambda_0 = -\d \mu_a \Lambda_0.
\end{equation}
Using $\Lambda_0$, various linear maps $\wh{W}_G\to B$ can be
described in terms of their integral kernels: For instance, the kernel
to the quantization map $\ca{Q}:\,\wh{W}_G\to \wh{\W}_G$ is the
element $\ca{Q}(\Lambda_0)\in \wh{\W}_G\otimes \Om(\g^*)$ and the
kernel of the restriction map $R_G^T:\,\,\wh{W}_G\to \wh{W}_T$ is the
element $R_G^T(\Lambda_0)$.

\subsection{Duality for $\wh{\W}_G$}\label{Sec:Duality2}
Just as for the commutative Weil algebra, the non-commutative Weil
differential can be interpreted in terms of a duality.

Let $\theta=\theta_a e_a$ be the 
left invariant and $\olt=\olt_a e_a$ the right invariant
Maurer-Cartan forms on $G$.  Denote by
$$ \eta=\f{1}{12}f_{abc}\theta_a\theta_b\theta_c $$
the canonical 3-form on $G$. It is closed, bi-invariant, and 
satisfies $\iota_a\eta=-\hh\d(\theta_a+\olt_a)$.

Suppose the group $G$ is a direct product of a connected, simply
connected Lie group and a torus.  
Then the canonical embedding 
$\g\to \mf{so}(\g)\subset \Cl(\g)$
given by the adjoint action exponentiates to a map 
$\tau:\,G\to \Spin(\g)\subset \Cl(\g)$. 
Explicitly, the map 
$\tau\in C^\infty(G)\otimes\Cl(\g)$ is given by
the formula (cf. \cite{al:no})
$$
\tau(\exp \mu) = \exp\left( - \frac{1}{2} f_{abc} \mu_a x_b x_c \right).
$$
The function $\tau$ defines an isomorphism of $\wh{\W}_G$ by
multiplication from the left. Using the symbol map, identify
$$ \wh{\W}_G\cong \ca{E}'(G)\otimes \wedge\g.$$
Using the forms $\theta_a$ we trivialize $T^*G\cong G\times \g^*$ 
and identify 
$$ \Om(G)\cong C^\infty(G)\otimes \wedge\g^*.$$
Let $\bfl A,\beta\bfr$ be the natural pairing given by these
identifications, and let  $\l A,\beta\r:=\bfl \tau^{-1} A,\beta\bfr$ be
the modified pairing.  
One finds (cf. \cite{al:no}, Proposition 5.6):
\begin{equation}\label{DualQWeil}
\d^\W=-(\d^{Rh}+\eta)^t,\ \ 
\iota_a=-(\iota_a^{Rh}+\hh(\theta_a+\olt_a))^t
\end{equation}
where $\d^{Rh},\iota_a^{Rh}$ are de Rham differentials and 
contractions on the space $\Om(G)$. 

The kernel of the identity map reads
\begin{equation}\label{Def:Lambda}
\Lambda:= e^{-x_a \bar{\theta}_a} \tau(g) \delta_g \in
\wh{\W}_G\otimes \Omega(G) ,
\end{equation}
and because of \eqref{DualQWeil} it has the properties,
\begin{equation}
\label{PropertiesLambda}
(\d^\W\otimes 1+1\otimes \d^{Rh}) \Lambda = \eta \Lambda , \ 
(\iota_a \otimes 1+1\otimes \iota_a^{Rh}) \Lambda = 
\frac{1}{2}(\theta_a + \bar{\theta}_a) \Lambda.
\end{equation}
As shown in Section 5 of 
\cite{al:no}, the 
product structure in $\wh{\W}_G$ can be expressed in 
terms of the following property of the form $\Lambda$:
\begin{equation}\label{MultLambda} 
\Lambda^1\Lambda^2=e^{-\hh \theta_a^1\olt_a^2}\Mult_G^*(\Lambda).
\end{equation}
Here $\Lambda^j\in\wh{\W}_G\otimes \Om(G\times G)$ are the 
pull-backs of $\Lambda$ to the $j$th $G$-factor, and  
$\Mult_G:\,G\times G\to G$ is group multiplication. 
The quantization map $\ca{Q}$ has 
a very simple description in terms of a dual map 
$\Om(G)\to \Om(\g)$: In fact, by Section 7 in \cite{al:no},
$$ \ca{Q}=(e^\varpi \circ \exp^*)^t $$
where $\varpi\in\Om^2(\g)^G$ is the image of $\eta$ under the 
de Rham homotopy operator $\Om^k(\g)\to \Om^{k-1}(\g)$.

\subsection{Duality for $\wh{\W}_T^-$}
View $\delta_t^-$ as a function of $t$ with values in $\ca{E}'(T)^-$.
As such it is the integral kernel of the identity map 
$\ca{E}'(T)^-\to \ca{E}'(T)^-$. The pairing 
between $\wh{\W}_T$ and $\Om(T)$ restricts to a pairing between 
$\wh{\W}_T^-$ and $\Om(T)$. If we denote by $\theta_k^T$ the 
Maurer-Cartan forms for $T$, the kernel of the identity map
$ \wh{\W}_T^-\to \wh{\W}_T^-$ with respect to this pairing is 
$$ \Lambda_T^-:=e^{-x_k\theta_k^T}\,\delta_t^-$$
One has the identity 
$$ (\Lambda_T^-)^1(\Lambda_T^-)^2=
e^{-\hh \theta_k^{T,1}\theta_k^{T,2}}
\Mult_T^* \Lambda_T^-$$
where $(\Lambda_T^-)^j\in \wh{\W}_T^-\otimes \Om(T\times T)$
are the pull-backs of $\Lambda_T^-$ to the $j$th factor, and 
$\Mult_T:T\times T\to T$ is the multiplication map.

\subsection{The integral kernel of the restriction map}
The statement of Theorem \ref{restriction} that $\ca{R}_G^T$ is a 
homomorphism of $T$-differential spaces 
can be rephrased in terms of the kernel 
$$ \ca{R}_G^T(\Lambda)\in (\wh{\W}_T^- \otimes \Om(G))^T.$$
of the restriction map:
$$ (\d+\eta)\ca{R}_G^T(\Lambda)=0,\ \ 
(\iota_k+\hh(\theta_k+\olt_k))\ca{R}_G^T(\Lambda)=0.$$
Here we write $\d=\d^\W\otimes 1+1\otimes \d^{Rh}$ etc. 
The pull-back of $\ca{R}_G^T(\Lambda)$ to the maximal torus has the 
following decription.
\begin{proposition}\label{RestrictionKernel}
Let $\iota_T:\,T\hra G$ be the inclusion of the maximal torus.
Then 
$$
\iota_T^*\ca{R}_G^T(\La)= \Lambda_T^-.
$$
\end{proposition}

\begin{proof}
We want to apply the operator $\ca{R}_G^T=r_G^T\otimes p_G^T\circ
\ca{T}$ to 
$$ \iota_T^*\La=
e^{-\sum_{k=1}^r x_k \theta_k^T}\,\tau(t)\,\delta_t.$$
The term $\exp(-\sum_{k=1}^r x_k \theta_k^T)$
lives in $\Cl(\t)_\C\otimes\Om(T)$ and hence commutes with both the
operator $\ca{T}=\exp(\hh r_{ab}\iota_a\iota_b)$ on $\Cl(\g)_\C$ 
and with the operator $p_G^T:\,\Cl(\g)_\C\to\Cl(\t)_\C$. Therefore,
$$\iota_T^* \ca{R}_G^T(\La)=
e^{-\sum_{k=1}^r x_k \theta_k^T}\,
\big(p_G^T\circ \ca{\T}\big)\tau(t)   r_G^T(\delta_t).
$$
To compute $\big(p_G^T\circ\ca{T}\big)\tau(t)$, 
choose $\mu\in\t$ with $t=\exp(\mu)$, and express $\tau(t)$ in the 
Cartan-Weyl basis $e_k,e_\alpha,e_{-\alpha}$. Note that since
$x_\alpha^2=x_{-\alpha}^2=0$ and $[x_\alpha,x_{-\alpha}]=1$,  
the element $x_\alpha x_{-\alpha}-x_{-\alpha} x_\alpha$ squares to 
$1$. Together with 
$\l\mu,[e_\alpha,e_{-\alpha}]\r =\tpi \l\mu,\alpha\r$, 
this shows that 
\beq \tau(\exp\mu)&=&\exp\big( i\pi \sum_{\alpha\in\mf{R}_+} 
\l\mu,\alpha\r 
(x_\alpha x_{-\alpha}-x_{-\alpha} x_\alpha)\big)\\
&=&
\prod_{\alpha\in\mf{R}_+}\big(\cos(\pi{\l\mu,\alpha\r})+i\sin(\pi 
\l\mu,\alpha\r) (x_\alpha x_{-\alpha}-x_{-\alpha}x_\alpha)\big).
\eeq
Applying the operator $p_G^T\circ \ca{T}=p_G^T\circ \prod_{\alpha\in\mf{R}_+}(1+\hh \iota_\alpha\iota_{-\alpha})$, we obtain 
$$ p_G^T\circ\ca{T}(\tau(\exp\mu))
=\prod_{\alpha\in\mf{R}_+} \exp(i\pi \l\mu,\alpha\r)
= \exp(\tpi\l\mu,\rho\r)=t^\rho.$$
Finally, $t^\rho$ and $r_G^T(\delta_t)$ combine by Example \ref{deltaminus},
and we find 
$$ \iota_T^*\ca{R}_G^T(\Lambda)
=e^{-\sum_{k=1}^r x_k \theta_k^T} \delta_t^-=\Lambda_T^-.$$
\end{proof}

\subsection{Quantization commutes with restriction}
The following result is a very useful computational tool,
since it replaces the complicated restriction map $\ca{R}_G^T$
and the non-abelian quantization map $\ca{Q}_G$ by the much simpler
restriction map $R_G^T$ and the abelian quantization map $\ca{Q}_T$.
\begin{theorem}\label{Th:QR=0}
The diagram 
\begin{equation}\label{BigDiagram}
\vcenter{
\xymatrix{
\wh{H}_G(M)\ar[d]^{R_G^T} \ar[r]^{\ca{Q}_G} & \wh{\H}_G(M)
\ar[d]^{\ca{R}_G^T}\\
\wh{H}_T(M)\ar[r]^{\pi^-\circ \ca{Q}_T} &  {\wh{\H}_T^-(M)}
}
}
\end{equation}
commutes.
\end{theorem}
\begin{proof}
The proof is similar to the proof in \cite{al:no}, Section 7,
that $\ca{Q}$ induces a ring homomorphism in cohomology. 
Consider the following two homomorphisms of $T$-differential spaces, 
$\phi_j:\,\wh{W}_G\to \wh{\W}_T^-$, 
$$ \phi_1= \pi^-\circ \ca{Q}_T \circ R_G^T,\ \ \ \ 
\phi_2= \ca{R}_G^T\circ \ca{Q}_G.$$
It suffices to show that $\phi_1,\phi_2$ are $T$-chain homotopic, that
is, there exists a $T$-equivariant, odd 
linear map $\psi:\,\wh{W}_G\to \wh{\W}_T^-$
such that $\psi$ commutes with Lie derivatives and contractions and
such that $\d^\W\circ \psi+\psi\circ \d^W=\phi_1-\phi_2$. 
To construct $\psi$ we
describe all maps in terms of their kernels.  The kernels of the maps
$\phi_j$ are
$$\ca{K}_j=(\phi_j \otimes 1)(\Lambda_0)\in (\wh{\W}_T^-\otimes
\Om(\g^*))^T.$$
Since $\phi_j$ are homomorphisms of $T$-differential spaces, they have
the property
$$  \d\ca{K}_j=0,\ \ \iota_k\ca{K}_j=-(\d\mu_k)\ca{K}_j.$$
The kernel $\ca{K}_1$ is invertible in the algebra 
$\wh{\W}_T^-\otimes \Om(\g^*)$, because its form degree $0$ part
$$(\ca{K}_1)_{[0]}|_\mu=\phi_1(\Lambda_0)_{[0]}|_\mu=
\delta_{\exp(\on{pr}_{\t^*}(\mu))}^-$$
is invertible at all points $\mu\in\g^*$. The element
$\ca{K}_1^{-1}\ca{K}_2$ then has properties
$$\iota_k(\ca{K}_1^{-1}\ca{K}_2)=0,\ \ \d(\ca{K}_1^{-1}\ca{K}_2)=0$$
so that $\ca{K}_1^{-1}\ca{K}_2$ is a closed element in
$(\wh{\W}_T^-\otimes \Om(\g^*))_{basic}$. Since $\g^*$ is
equivariantly contractible, the pull-back under the inclusion
$\iota:\,\{0\}\to \g^*$ induces an isomorphism in cohomology.  Since
both $\phi_j$ map the identity to the identity, $\iota^*\ca{K}_j=1$
and therefore $\iota^*(\ca{K}_1^{-1}\ca{K}_2)=1$. Hence, there exists
$\ca{N}'\in (\wh{\W}_T^-\otimes \Om(\g^*))_{basic}$ with
$$\ca{K}_1^{-1}\ca{K}_2=1-\d\ca{N}'.$$
and letting $\ca{N}=\ca{K}_1\ca{N}'$ and using that $\ca{K}_1$ is
closed we find
$$ \ca{K}_2=\ca{K}_1-\d \ca{N}.$$
Consequently, the linear map $\psi$ with kernel $\ca{N}$ is
$T$-equivariant and satisfies 
$\d^\W\circ \psi+\psi\circ \d^W=\phi_1-\phi_2$.  Since
$$ (\iota_k+\d\mu_k)\ca{N}=0$$
the map $\psi$ commutes with $T$-contractions.
\end{proof}

Taking $M=\pt$, one recovers the fact from Lie group theory that the
diagram
\begin{equation}\label{quantcom}
\vcenter{\xymatrix{ \ca{E}'(\g)^G \ar[r]^{\on{Duf}} \ar[d]^{r_G^T} &
\ca{E}'(G)^G \ar[d]^{r_G^T}\\ \ca{E}'(\t) \ar[r]^{\pi^-\circ \exp_*} &
{\ca{E}'(T)^-}}}.
\end{equation} 
commutes. 

\subsection{Ring structure}\label{RingStructure}
It is known that all arrows in the diagram \eqref{quantcom} are in
fact ring maps: For the horizontal arrows this is the Duflo theorem,
while for the vertical arrows it follows easily from the formula for
the convolution of irreducible characters.  This raises the question
whether the arrows in the diagram \eqref{BigDiagram} are ring maps as
well? The left vertical arrow $R_G^T:\,\wh{H}_G(M)\to \wh{H}_T(M)$ is
a ring map, and the horizontal arrows are ring maps by Theorem 8.1 in
\cite{al:no}.

To address the problem whether the right vertical arrow 
$\ca{R}_G^T:\,\wh{\H}_G(M)\to
\wh{\H}_T^-(M)$ is a ring map, let us consider the two maps
$\phi_j:\,\wh{\W}_G\otimes \wh{\W}_G\to \wh{\W}_T^-$ defined
respectively by restriction after multiplication and multiplication
after restriction. The kernels $\ca{K}_j\in \wh{\W}_T^-\otimes
\Om(G\times G)$ for these two maps are
$$ \ca{K}_1=\ca{R}_G^T(\Lambda^1\Lambda^2),\ \ \ 
\ca{K}_2=\ca{R}_G^T(\Lambda^1)\ca{R}_G^T(\Lambda^2)$$
where the superscripts denote the pull-backs to the 
respective $G$-factor. Since $\phi_j$ are $T$-equivariant chain 
maps commuting with $T$-contractions, 
these kernels satisfy 
$$ (\d+\eta^1+\eta^2)\ca{K}_j=0,\ \ 
(\iota_k+\hh(\theta_k^1+\olt_k^1)
+\hh(\theta_k^2+\olt_k^2))\ca{K}_j=0
$$ 
($k=1,\ldots,r$). To show that $\phi_1,\phi_2$ induce the same map in
cohomology, it would be enough to show that there exists $\ca{N}\in
(\wh{W}_T^-\otimes \Om(G\times G))^T$ with
\begin{equation}\label{A}
 \ca{K}_2=\ca{K}_1-(\d+\eta^1+\eta^2)\ca{N}
\end{equation}
and 
\begin{equation}\label{B}
(\iota_k+\hh(\theta_k^1+\olt_k^1)+\hh(\theta_k^2+\olt_k^2))\ca{N}=0.
\end{equation}
We do not know whether $\ca{N}$ with these properties 
exists in general, however, for our applications the 
following weaker statement will suffice:
\begin{proposition}
There exists a $T$-invariant tubular neighborhood $V$ of $T$ in $G$, and 
$\ca{N}\in (\wh{\W}_T^-\otimes \Om(V\times V))^T$ such that 
equations \eqref{A} and \eqref{B} hold over $V\times V$.
\end{proposition}
\begin{proof}
We begin by showing that the pull-backs of $\ca{K}_1$ and $\ca{K}_2$
under the inclusion $\iota_{T\times T}:\,T\times T \to G\times G$
coincide. 
We compute 
$\iota_{T\times T}^*\ca{K}_2$ using 
Proposition \ref{RestrictionKernel}:
$$
\iota_{T\times T}^*\ca{K}_2=
(\Lambda_T^-)^1(\Lambda_T^-)^2
=e^{-\hh \sum_{k=1}^r \theta_k^{1,T} \theta_k^{2,T}}
\Mult_T^*\Lambda_T^-
$$
By another 
application of Proposition \ref{RestrictionKernel}
we find 
$$
\iota_{T\times T}^*\ca{K}_1
=\iota_{T\times T}^* e^{-\hh \theta_a^1\olt_a^2} \Mult_G^*\ca{R}_G^T(\Lambda)
=
e^{-\hh \sum_{k=1}^r
\theta_k^{1,T} \theta_k^{2,T}} \Mult_T^*\Lambda_T^-$$
Since $\iota_{T\times T}^*\ca{K}_1$ has non-vanishing form degree 
$0$ part, there exists a $T$-invariant neighborhood $V$ of $T$ 
such that $\ca{K}_1$ has non-vanishing form degree $0$ part, 
and is therefore invertible, on $V\times V$.  
We then have, over $V\times V$, 
$$ \d(\ca{K}_1^{-1}\ca{K}_2)=0,\ \ \iota_k(\ca{K}_1^{-1}\ca{K}_2)=0
$$
while $\iota_{T\times T}^*(\ca{K}_1^{-1}\ca{K}_2)=1$. 
Since $V$ retracts $T$-equivariantly onto $T$, the de Rham 
homotopy operator for the retraction $V\times V\to T\times T$ 
defines an element  
$\ca{N}'\in (\wh{\W}_T^-\otimes \Om(V\times V))_{basic}$ 
such that 
$$ \ca{K}_1^{-1}\,\ca{K}_2=1-\d\ca{N}'.$$
Letting $\ca{N}=\ca{K}_1\ca{N}'$ this gives \eqref{A}. 
Equation \eqref{B} follows from the equation for 
$\ca{K}_2,\ca{K}_1$. 
\end{proof}

We will apply this result in the following way. 
Suppose $\beta_1,\beta_2\in (\wh{\W}_G\otimes \Om(M))_{basic}$
are two equivariant cocycles for a given $G$-manifold $M$. 
Suppose moreover that for some $T$-invariant open subset 
$U\subset M$, the restrictions $\beta_j|U$ take values in 
$(\ca{E}'(V)\otimes \Cl(\g)\otimes \Om(U))$ where $V$ is 
the tubular neighborhood of $T$ as in the proposition. 
Then the two cocycles
$\ca{R}_G^T(\beta_1\beta_2)$ and
$\ca{R}_G^T(\beta_1)\ca{R}_G^T(\beta_2)$ in the space 
$(\wh{\W}_T^-\otimes \Om(M))_{basic}$ 
are cohomologous if restricted to $U$.

\section{Localization formulas for $\wh{\H}_G(M)$.}
\label{Sec:LocForm}
In this Section we derive a new nonabelian localization formula which
applies to equivariant cohomology $\wh{\H}_G(M)$.
\subsection{Localization formulas}
The localization formulas discussed in this paper are derived from a
common principle, Theorem \ref{th:loc} below.  Let $M$ be an oriented
manifold, acted upon by a compact Lie group $G$. For any $\xi\in \g$
the derivation
$$ \d_\xi=\d-\tpi\iota_\xi:\,\Om(M)\to\Om(M)$$
is a differential on the kernel of $L_\xi$:
$$C_\xi(M)=\{\beta\in\Om(M),\,L_\xi\beta=0\}.$$ 
There is a natural ring homomorphism
$$ J_\xi:\ C_G(M)\to C_\xi(M),\,\beta\mapsto \beta(\tpi\xi). $$
Here we view elements of $C_G(M)$ as polynomials on $\g_\C$ with
values in differential forms, and $\beta(\tpi\xi)$ is the evaluation
of that polynomial.  The localization formula expresses the integral
of $\d_\xi$-cocycles in terms of integrals over the zeroes of the
vector field $\xi_M$.  Since $G$ is compact, $\xi_M^{-1}(0)$ is a
union of compact, smooth, orientable submanifolds. Let $\F(\xi)$
denote the set of connected components of $\xi_M^{-1}(0)$, and for any
$F\in \F(\xi)$ choose an orientation on $F$. Let the normal bundle
$\nu_F$ be equipped with the orientation obtained from the orientation
on $M$ and on $F$.

Let $T$ be a maximal torus of $G$ with $\xi\in\t$. Then every $F$ is a
$T$-invariant submanifold; let $\Eul(\nu_F)\in C_T(F)$ be the
$T$-equivariant Euler form of $\nu_F$, for some choice of invariant
metric and connection. We claim that
$$J_\xi(\Eul(\nu_F))=\Eul(\nu_F,\tpi \xi)\in C_\xi(F) $$
is {\it invertible} in $C_\xi(F)$.  For this it suffices
to note that its form degree $0$ part is invertible. 
Up to a constant, the form degree $0$ part
is the Pfaffian of the endomorphism of the fibers of
$\nu_F$ defined by $\xi$. Since $\xi_{\nu_F}$ vanishes 
precisely on $F$, this endomorphism is
invertible, and therefore the Pfaffian is non-zero.
 
The following integration formula \eqref{bvloc} was proved by
Berline-Vergne \cite{be:ze} in the Cartan model of equivariant
cohomology and by Atiyah-Bott \cite{at:mom} using the topological
definition.  It generalizes earlier work of Bott \cite{bo:ve,bo:re} on
the zeroes of holomorphic vector fields on complex manifolds and of
Baum-Cheeger \cite{ba:in} (see also \cite{ko:tr}, Chapter II.6) on
zeroes of Killing vector fields on Riemannian manifolds.

We include the proof, since the result in 
\cite{be:ze} or \cite{at:mom}
is stated in less generality than what is needed
here, although their argument extends. 

\begin{theorem}[Localization formula]\label{th:loc}
Let $M$ be an oriented $G$-manifold, and $\xi\in\g$.  For any
compactly supported $\d_\xi$-cocycle $\gamma\in C_\xi(M)$, the
integral of $\gamma$ over $M$ is given by the localization formula,
\begin{equation} \label{bvloc}
\int_M \gamma = \sum_{F \in  \F(\xi)} \ \int_F
\frac{\iota_F^* \gamma}{\Eul(\nu_F, 2\pi i \xi) }.
\end{equation}
\end{theorem}
\begin{proof}
Our proof is a simple adaptation of the proof given in 
Guillemin-Sternberg \cite{gu:su}, which in turn is based 
on the original proofs of Berline-Vergne and Atiyah-Bott.
We may assume that $G=T$ is abelian.  

For each $F$, let $U_F\cong \nu_F$ 
be a $T$-invariant tubular neighborhood of
$F$. Let $f\in C^\infty(M,\R)$ be a function such that 
$f$ vanishes in a neighborhood of the fixed point set 
and $f=1$ in a neighborhood of $M\backslash \bigcup_F U_F$. 
Choose an invariant metric $g$ on $M$ and let 
$\psi$ be a $T$-invariant 1-form, defined on the complement of the 
fixed point set, $\psi:=g(\xi_M,\cdot)/g(\xi_M,\xi_M)$. 
Then $\d_\xi\psi$ is invertible. Write 
$$\gamma=\gamma'+\d_\xi\big(f\f{\psi}{\d_\xi\psi}
\gamma\big).$$
Then $\int_M \gamma=\int_M\gamma'$. The form $\gamma'$ 
is supported on $\bigcup_F U_F$ and agrees with $\gamma$ 
near $\bigcup_F F$. 

This reduces the proof to the case where $\gamma$ is compactly
supported in $M=\nu_F$ for some $F$. Let $\on{Th}(\nu_F)\in
C_T(\nu_F)$ be the Mathai-Quillen form \cite{ma:th} 
representing the T-equivariant Thom 
class of $F$. It has compact support in fiber direction and its
pull-back to $F$ represents the equivariant Euler form,
$\Eul(\nu_F)$. Let $\on{Th}(\nu_F,\tpi\xi)=J_\xi \on{Th}(\nu_F)$.
Then the pull-back of $\Eul(\nu_F,\tpi\xi)$ to $\pi_F:\nu_F\to F$ is
$\d_\xi$-cohomologous to $\on{Th}(\nu_F,\tpi\xi)$. Consequently the
quotient $\on{Th}(\nu_F,\tpi\xi)/\pi_F^*\Eul(\nu_F,\tpi\xi)$ is
$\d_\xi$-cohomologous to $1$. Using the main property of the Thom form
that its fiber integral is $1$, and replacing $\gamma$ by the
$\d_\xi$-cohomologous form $\pi_F^*\iota_F^*\gamma$ we
find:
\beq 
\int_{\nu_F}\gamma &=&\int_{\nu_F}\f{\gamma}{\pi_F^*\Eul(\nu_F,\tpi\xi)}
\on{Th}(\nu_F,\tpi\xi)\\&=&
\int_{\nu_F}
\pi_F^* \big(\f{\iota_F^* \gamma}{\Eul(\nu_F,\tpi\xi)}\big)
\on{Th}(\nu_F,\tpi\xi)
=\int_F \f{\iota_F^*\gamma}{\Eul(\nu_F,\tpi\xi)}.\eeq
\end{proof}

\begin{remark}\label{wedge}
Later we need the following observation about the right hand side
of the integration formula. Suppose $\beta\in C_\xi(F)$ 
is a compactly supported $T$-invariant $\d_\xi$-cocycle, and 
$\zeta\in\t$. Then $\iota_\zeta\beta$ is a $\d_\xi$-cocycle: 
$\d_\xi\iota_\zeta\beta=-\iota_\zeta\d_\xi\beta+L_\zeta\beta=0$. 
%
%
We claim that
$$ \int_F \f{\iota_\zeta \beta}{\Eul(\nu_F,\tpi\xi)}=0.$$
As in the proof, $\pi_F^*\beta\in C_\xi(\nu_F)$ can be 
replaced with a compactly supported, $T$-invariant, 
$\d_\xi$-cohomologous form $\gamma$, which agrees with 
$\pi_F^*\beta$ near $F$. By the localization formula, used 
in reverse, the above integral equals $\int_{\nu_F} 
\iota_\zeta\gamma$, which is $0$ since $\iota_\zeta\gamma$
has vanishing top degree.  
\end{remark}

The localization formula is usually stated as an integration 
formula for cocycles in  $C_G(M)$, 
taking $\gamma$ of the form $\beta(\tpi\xi)$ 
for some $\beta\in C_G(M)$. 
However, not every $\d_\xi$-cocycle arises in this way, 
so \eqref{bvloc} is a more general statement.
Indeed, the complexes $\wh{C}_G(M)$ and $\wh{\ca{C}}_G(M)$ produce 
other examples:

\subsubsection{The map $J_\xi:\ \wh{C}_G(M)\to C_\xi(M)$.}
\label{11}
For any $\xi\in\g$, the map 
$$J_\xi:\  \wh{C}_G(M)\to C_\xi(M),\ \beta\mapsto 
\l\beta,e^{\tpi\l\mu,\xi\r}\r$$
is a chain map since
\beq 
(\d_M-\tpi\iota(\xi_M))\l\beta,e^{\tpi\l\mu,\xi\r}\r
&=&\l (v^a-\tpi \xi^a)\iota_a\beta,e^{\tpi\l\mu,\xi\r}\r
\\
&=&\l \iota_a\beta,(\f{\p}{\p\mu_a}-\tpi \xi^a)e^{\tpi\l\mu,\xi\r}\r
=0.
\eeq
The localization formula describes the Fourier
components of the integration map 
$$\int_M:\,\wh{H}_G(M)\to \wh{H}_G(\pt)=\ca{E}'(\g^*)^G .$$

\subsubsection{The map ${\ca{J}}_\lambda:\ \wh{\ca{C}}_T(M)\to 
C_{\lambda}(M)$.}
\label{22}
Suppose $G=T$ is a torus, and that its Lie algebra $\t$ is equipped
with an inner product.  For any weight $\lambda\in\Lambda^*\subset\t$
the map
$$ \ca{J}_\lambda:\ \wh{\ca{C}}_T(M)\to C_\lambda(M),\ \ 
\beta\mapsto \l\beta,t^\lambda\r$$
is a chain map, 
by a calculation similar to that for $\wh{\ca{C}}_G(M)$:
\beq 
(\d_M-\tpi\iota(\lambda_M))\l\beta,t^\lambda\r
&=&\l (v^a-\tpi \lambda^a)\iota_a\beta,t^\lambda\r \\
&=&\l \iota_a\beta,(e_a^L-\tpi \lambda_a)t^\lambda\r=0.
\eeq
We hence obtain a 
map $J_\lambda:\ \wh{\ca{C}}_T(M)\to C_{\lambda}(M)$.
In this case, the localization formula 
describes the Fourier coefficients $\l\int_M\beta,t^\lambda\r$
of the integration map 
$$\int_M:\wh{\ca{C}}_T(M)\to \wh{\ca{C}}_T(\pt)=\ca{E}'(T).$$

\subsubsection{The map $\ca{J}_\lambda:\ 
\wh{\ca{C}}_G(M)\to C_{\lambda+\rho}(M)$.}
The map from \ref{22} generalizes to the non-abelian case as 
follows. We assume that $G$ is a product of a simply connected 
group and a torus. The main result of Section \ref{Sec:Abel}
was the construction of a chain map, $\wh{\ca{C}}_G(M) \to 
\wh{\ca{C}}_T^-(M)$. Composing with the inclusion 
$\wh{\ca{C}}_T^-(M)\to \wh{\ca{C}}_T(M)$ and the map 
$\wh{\ca{C}}_T(M)\to C_{\lambda+\rho}(M)$ we obtain a chain map, 
\begin{equation} 
{\J}_\lambda:\,\wh{\ca{C}}_G(M)\to C_{\lambda+\rho}(M),
\beta\mapsto  e^{\hh \iota(r_M)}\l \beta,\Delta_\lambda\r.
\end{equation}
Since the contraction operator $e^{\hh \iota(r_M)}$ does not change the 
top form degree part, 
$$ \int_M \ca{J}_\lambda(\beta)
=\l \int_M \beta,\Delta_\lambda\r
=(\dim V_\lambda)^{-1} \l \int_M\beta,\chi_\lambda\r 
$$
showing that the localization formula describes the Fourier
coefficients of the map
$$ \int_M:\ \wh{\H}_G(M)\to \wh{\H}_G(\pt)=\ca{E}'(G)^G.$$

Applying Theorem  \ref{th:loc} to this situation we have proved 
the main result of this paper:
\begin{theorem}\label{Th:main}
Let $\beta\in\wh{\ca{C}}_G(M)$ be an equivariant cocycle. 
Then 
$$ \l \int_M\beta,\chi_\lambda\r=
\dim V_\lambda\sum_{F\in\F(\lambda+\rho)}
\int_F \f{\iota_F^*(e^{\hh \iota(r_M)} \l\beta,\S_\lambda\r)}
{\Eul(\nu_F,\tpi(\lambda+\rho))}.
$$
for all $\lambda\in\Lambda^*_+$. 
\end{theorem}
For the subcomplex $\ca{C}_G(M)=(U(\g)_\C \otimes \Om(M))^G$, 
the integration map takes values 
in the space of Casimir elements $U(\g)_\C ^G$. Any Casimir 
element $P$ is determined by its eigenvalues in all irreducible 
representations, or equivalently by its traces 
$\on{Tr}_{V_\lambda}(P)$.

%
%
%

Below we will need the following additional facts about the 
map $\ca{J}_\lambda$. First, 
suppose $\beta_0\in\wh{C}_G(M)$
is a cocycle and that
$\beta=\ca{Q}(\beta_0)$ is its quantization. Then  
then ${\J}_\lambda(\beta)$ and $J_{\lambda+\rho}(\beta_0)$
are cohomologous in $C_{\lambda+\rho}$.  
This follows from the ``quantization commutes with restriction'' 
Theorem, \ref{Th:QR=0} and the definition of the maps
$J_{\lambda+\rho}$ and ${\J}_{\lambda}$. 
Theorem \ref{Th:main} does not yield anything new for such 
classes. 

Next, suppose $\beta\in\wh{\ca{C}}_G(M)$ 
is a product of two cocycles, $\beta=\beta_1\odot\beta_2$. 
Suppose that both $\beta_1$ and $\beta_2$ have the property that 
if the support of the 
test function $\phi\in C^\infty(G)$ does not meet 
$T$, then $\l\beta_j,\phi\r$ vanish near all $F\in \F(\lambda+\rho)$. 
(For example, this is trivially the case for any 
$\beta_j$ contained in $\ca{C}_G(M)\subset \wh{\ca{C}}_G(M)$.
According to the remarks at the end of Section \ref{RingStructure}, 
this means that 
$$(r_G^T\otimes e^{\hh\iota(r_M)})(\beta_1\odot\beta_2)\ \ \ 
\mbox{ and }\ \ \ 
(r_G^T\otimes e^{\hh\iota(r_M)})(\beta_1)\odot
(r_G^T\otimes e^{\hh\iota(r_M)})
(\beta_2)
$$
are cohomologous near the fixed point set. 
Here we are using the product $\odot$ on $\wh{\ca{C}}_T(M)^-$
the product  induced by the multiplication 
$$ \diag_M^*\circ \exp(-\hh \sum_{k=1}^r \iota_k^1\iota_k^2)$$
on $\Om(M)^T$. It follows that in the localization formula, 
one can replace 
$\J_\lambda(\beta_1\odot\beta_2)$ by $\J_\lambda(\beta_1)\odot 
\J_\lambda(\beta_2)$. 
According to Remark \ref{wedge}, in the integral 
over $F$ we may replace $\odot$ by simply the wedge product. 
Hence we have shown: 
\begin{lemma}\label{Lem:Prod}
Let $M$ be a compact oriented $G$-manifold, and 
$\lambda\in\Lambda^*_+$. Suppose the cocycles 
$\beta_1,\beta_2\in \wh{\ca{C}}_G(M)$ have the property that 
if the support of $\phi\in C^\infty(G)$ does not meet $T$, 
then the support of $\l\beta_j,\phi\r$ does not meet the
fixed point sets $F\in\F(\lambda+\rho)$. Then 
$$ \l\int_M(\beta_1\odot\beta_2),\chi_\lambda\r
=\dim V_\lambda \sum_{F\in\F(\lambda+\rho)}\int_F 
\f{\J_\lambda(\beta_1)\J_\lambda(\beta_2)}{\Eul(\nu_F,\tpi(\lambda+\rho))}.
$$
\end{lemma}

\section{Applications to moment maps}
\label{Sec:Appl}
In this Section we discuss applications of our localization
formula to moment map theory. The main result is a 
Duistermaat-Heckman formula for group-valued moment maps.

\subsection{Duistermaat-Heckman formula for $\g^*$-valued moment maps}
We recall the definition of a Hamiltonian $G$-space.
\begin{definition}
A Hamiltonian $G$-space is a triple $(M, \omega_0, \Phi_0)$ consisting
of a $G$-space $M$, a 2-form $\omega_0$, and an equivariant map
$\Phi_0 \in C^\infty(M, \g^*)^G$ satisfying
\begin{equation} \label{HGspace}
\begin{array}{cl}
\d \omega_0 =0 & \on{Cocycle \ condition} \\
\iota(\xi_M) \omega_0 = \d \l \Phi_0 , \xi \r & 
\on{Moment \ map \ condition} \\
\ker((\omega_0)_x) = \{ 0 \} & 
\on{ Non-degeneracy}
\end{array}
\end{equation}
\end{definition}
The first two conditions imply that the 
form $\omega_0$ is $G$-invariant.
Basic examples of Hamiltonian $G$-spaces are coadjoint orbits
$\O=G\cdot\mu$ for $\mu\in\g^*$, with moment map the
embedding $\Phi_0: \O \to \g^*$.

Given a Hamiltonian $G$-space $(M, \omega_0, \Phi_0)$ and a 
coadjoint orbit $\O=G\cdot\mu$ such that $\mu$ is 
a regular value of the moment map 
the reduced space $M_\O := \Phi^{-1}(\O)/G$
is an orbifold with a 
naturally induced symplectic form $(\om_0)_\O$.

For any Hamiltonian $G$-space $M$ let the 
equivariant Liouville form be the element of 
$\wh{W}_G \otimes \Omega(M)$ defined as 
\begin{equation}
\L_0 = e^{\omega_0} \Phi_0^* \Lambda_0, 
\end{equation}
where $\Lambda_0\in \wh{\W}_G\otimes\Om(\g^*)$ 
is the integral kernel introduced in \eqref{Def:Lambda0}. 
Equations \eqref{HGspace} and \eqref{Lambda0} show that 
$\L_0$ is a cocycle contained in the basic subcomplex 
$(\wh{W}_G\otimes \Om(M))_{basic}$.
The corresponding form $\L_0^{Car}=P_{hor}^W(\L_0)$ in the Cartan model is 
$$ 
\L_0^{Car}=e^{\om_0} \delta_{\Phi_0}\in \wh{{C}}_G(M).$$
The top degree component of its integral over $\g^*$,
$$
\Gamma_0:= \l  \L_0^{Car}, 1 \r_{[top]}= \left( e^{\om_0} \right)_{[top]}
$$
is a volume form on $M$, defining in particular a natural orientation.
For a compact oriented Hamiltonian $G$-space $(M, \om_0, \Phi_0)$
with Liouville form $\L_0$ one defines the Duistermaat-Heckman measure
(DH distribution) as the integral,
\begin{equation}
\mathfrak{m}_0:=\int_M \L_0^{Car}
\in  \ca{E}'(\g^*)^G.
\end{equation}
The support of the distribution $\mf{m}_0$ is contained in 
the image of $\Phi_0$ and its singular support in the set of singular
values. These two properties extend to more general 
DH-distributions constructed from equivariant cocycles 
$\beta_0 \in C_G(M)$
\begin{equation}
\mathfrak{m}^{\beta_0}_0 := \int_M \beta_0 \L_0^{Car} 
\in \ca{E}'(\g^*)^G.
\end{equation}
The main motivation for studying the twisted DH distributions
$\m_0^{\beta_0}$ is that they encode intersection pairings on reduced
spaces.  
More precisely, Let $\Vol_{\g^*}$ be the Riemannian volume volume 
on $\g^*$ for a given invariant inner product on $\g$, and 
let $\Vol G$ be the corresponding Riemannian volume of the group $G$.

For any coadjoint orbit $\O$ in the set of
regular values of $\Phi_0$, the value of this distribution at $\O$ is
given by the formula
\begin{equation}  \label{pairings}
 \frac{\mathfrak{m}^{\beta_0}_0}{\d \Vol_\g^*} \Big|_\O =
\frac{\Vol G}{k \Vol \O } \int_{M_\O} \kappa_\O(\beta_0) e^{(\om_0)_\O}.
\end{equation}
Here $k$ is the cardinality of a generic stabilizer for the $G$-action
on $M$, $\Vol \O$ is the Liouville volume of the coadjoint orbit and
$\kappa_\O: C_G(M) \to C(M_\O)$ is the chain map given as a
composition of the pull-back map $C_G(M)\to C_G(\Phi_0^{-1}(\O))$ and
the Cartan map $C_G(\Phi_0^{-1}(\O))\to \Om(M_\O)$.  By a theorem of
Kirwan 
\cite{ki:coh}, the map $\kappa_\O$ induces a surjective map in
cohomology.

Applying the localization formula to $\beta_0\L_0^{Car}$, and 
using
$$J_\xi(\beta_0\L_0^{Car})=J_\xi(\beta_0) J_\xi(\L_0^{Car})
=\beta_0(\tpi\xi)
e^{\om_0+\tpi\l\Phi_0,\xi\r}$$
one obtains the following formula for the Fourier coefficients 
of $\mathfrak{m}^\beta_0$:
\begin{proposition}[Duistermaat-Heckman formula]
Let $G$ be a compact, connected Lie group, 
$(M, \om_0, \Phi_0)$ a compact Hamiltonian
$G$-manifold, and let $\beta_0 \in C_G(M)$ be a cocycle.
Then the Fourier components of the twisted DH distribution 
$\mathfrak{m}^{\beta_0}_0$ are given by the formula,
\begin{equation}\label{DH1}
\l \mathfrak{m}^{\beta_0}_0, e^{2\pi i \l \mu, \xi \r} \r =
\sum_{F \in \F(\xi)}  e^{\tpi\l \Phi_{0,F}, \xi \r}
\int_F \frac{i_F^* \beta_0(2\pi i\xi) e^{\omega_{0,F}}}
{\Eul(\nu_F, 2\pi i \xi)}  .
\end{equation}
Here $\om_{0,F}=\iota_F^*\om_0$ and $\Phi_{0,F}=\iota_F^*\Phi_0$. 
\end{proposition}
\subsection{Group valued moment maps}
In this section we prove an analogue of the Duistermaat-Heckman
formula for group valued moment maps.
\begin{definition}\cite{al:mom}
A $G$-valued Hamiltonian $G$-space is a triple $(M, \omega, \Phi)$
consisting of a $G$-space $M$, an invariant 
2-form $\omega$, and an equivariant
map $\Phi \in C^\infty(M, G)^G$ satisfying the following
conditions,
\begin{equation}
\begin{array}{cl}
\d \omega = \Phi^* \eta & \Phi-\on{relative\  cocycle \ condition} \\
\iota_a \omega = \frac{1}{2} \Phi^*(\theta_a + \bar{\theta}_a) & 
\on{Moment \ map \ condition} \\
\ker((\omega_0)_x) = \{ \xi_M(x)| \Ad_{\Phi(x)} \xi = - \xi \} & 
\on{ Minimal \ degeneracy \ condition }
\end{array}
\end{equation}

\end{definition}
For the motivation of this definition we refer to \cite{al:mom}.
Basic examples of Hamiltonian $G$-spaces are conjugacy classes
$\Co:= G\cdot g$ for $g\in G$, with moment map the
embedding $\Phi:\Co \to G$ .

Many concepts from the theory of Hamiltonian $G$-spaces have 
analogues in the group-valued setting. In particular, if 
$\Co\subset G$ is a conjugacy class contained in the set 
of regular values of $\Phi$, the 
reduced space $M_\Co:= \Phi^{-1}(\Co)/G$, is an 
orbifold with a naturally induced symplectic form  $\om_\Co$. 
Our aim is to extract information about intersection pairings on
reduced spaces from localization formulas on $M$.

We assume that the group $G$ is a direct product of a connected, 
simply connected Lie group and a torus (this assumption will be lifted in 
Section \ref{Sec:General}). The equivariant Liouville form 
associated to a $G$-valued
Hamiltonian $G$-space $(M, \omega, \Phi)$ is defined by
\begin{equation} \label{Lio}
\L:= e^{\omega}  \Phi^* \Lambda \in \wh{\W}_G \otimes \Omega(M),
\end{equation}
with $\Lambda$ as in \eqref{Def:Lambda}. 
The properties \eqref{PropertiesLambda}
of $\La$ imply that $\L$ is an equivariant 
cocycle in the basic subcomplex. 
The Liouville form in the Cartan model, 
$\L^{Car}=P_{hor}^W(\L)$, is given by a fairly complicated 
expression worked out in \cite{al:du}. It is also shown in 
\cite{al:du} that, similar to the 
$\g^*$-valued setting, the top degree part  
$\Gamma=\l \L^{Car},1\r_{[top]}$  
is a volume form on $M$, giving rise to a natural orientation. 

Again, one can define the DH distribution $\mathfrak{m}$ as an 
integral over $M$ of the equivariant Liouville form
$$
\mathfrak{m}:=\int_M \L^{Car} \in \ca{E}'(G)^G.
$$
For any cocycle 
$\beta\in\C_G(M)$, define more general twisted DH-distributions 
$$
\mathfrak{m}^\beta:= \int_M \beta \odot \L^{Car} \in \ca{E}'(G)^G.
$$
By Proposition 5.7 of \cite{al:du} $\m^\beta$ is supported on
the image of $\Phi$, and its singular support is contained in
the set of singular values of $\Phi$. The relationship between 
intersection pairings on the reduced space $M_\Co$, 
for $\Co$ contained in the set of regular values, and
twisted DH distributions $\mathfrak{m}^\beta$ is as follows,
(cf. \cite{al:du}):
\begin{equation}\label{interpairings}
\frac{\mathfrak{m}^\beta}{\d \Vol_G} \Big|_\Co =
\frac{\Vol G}{ k \Vol \mathcal{C}} \ 
\int_{M_\Co} \kappa_\Co(\beta) e^{\omega_\Co}.
\end{equation}
Here $\d \Vol_G$ is the Haar measure on $G$ defined by the inner
product on $\g$, $\Vol \mathcal{C}$ is the Liouville 
volume of the conjugacy
class $\Co$ with respect to the volume form $\Gamma_\Co$,  
and $k$ the number of elements
in a generic stabilizer. 
The map $\kappa_\Co:\,\ca{C}_G(M) \to \Om(M_\Co)$ is 
the composition of pull-back to $\Phinv(\Co)$, the 
isomorphism $\ca{C}_G(\Phinv(\Co))\cong C_G(\Phinv(\Co))$ given 
by the inverse of the quantization map, and the Cartan map
$C_G(\Phinv(\Co))\to M_\Co$. 

\subsection{DH formula for group valued moment maps}\label{DHGroup}
Similar to the case of Hamiltonian $G$-spaces, the Fourier
decomposition of the DH-distributions $\m^\beta$ can be obtained from 
a localization formula. To state the result, note first that for 
any $\lambda\in\Lambda^*_+$, all $F\in \F(\lambda+\rho)$ are contained
in the pre-image of the maximal torus $\Phinv(T)$, since 
$\Phi$ is equivariant and since the vector field for the conjugation 
action $(\lambda+\rho)_G$ vanishes exactly on $T$. 
Let $\om_F=\iota_F^*\om$ and $\Phi_F=\iota_F^*\Phi$. Since 
$\Phi_F$ takes values in $T$, we can compose with the map 
$T\to S^1,\,t\mapsto t^{\lambda+\rho}$ and this composition 
is denoted $(\Phi_F)^{\lambda+\rho}$. Indeed, $F$ is a 
group valued Hamiltonian $T$-space with $\om_F$
as 2-form and $\Phi_F$ as moment map. From the fact that 
$(\lambda+\rho)_M$ vanishes on $F$, one deduces that 
$(\Phi_F)^{\lambda+\rho}$ is constant along $F$.

\begin{theorem}[DH formula for group-valued moment maps] \label{qDH}
Let $G$ be a direct product of a connected, simply-connected Lie
group and a torus, $(M, \omega, \Phi)$ a compact $G$-valued
Hamiltonian $G$-space, and $\beta=\ca{Q}(\beta_0) 
\in C_G(M)$ be an equivariant cocycle. The
Fourier coefficients of the twisted DH distribution
$\m^\beta$ are given by the formula,
\begin{equation} \label{fqdh}
\l \m^\beta, \chi_\lambda \r=
\dim V_\lambda \sum_{F \in \F(\lambda+ \rho)} (\Phi_F)^{\lambda+\rho}
\int_F  \frac{\iota_F^* \beta_0(2\pi i(\lambda+\rho))}
{\Eul(\nu_F, 2\pi i(\lambda+\rho))}\ 
e^{\omega_F} .
\end{equation}
\end{theorem}

\begin{proof}
We apply the localization formula to the class
$\beta\odot \L^{Car}$. Note that the cocycles 
$\beta$ and $\L^{Car}$ satisfy the condition from 
Lemma \ref{Lem:Prod}, which therefore gives 
$$
\l\m^\beta,\chi_\lambda\r=\dim V_\lambda 
\sum_{F\in\ca{F}(\lambda+\rho)}\int_F
\f{    \ca{J}_{\lambda}(\beta) 
\ca{J}_{\lambda}(\L^{Car})  }
{\Eul(\nu_F,\tpi(\lambda+\rho))}.
$$
By the discussion preceding Lemma \ref{Lem:Prod}, 
$\ca{J}_\lambda(\beta)$ may be replaced with 
$$
J_{\lambda+\rho}(\beta_0)=
\beta_0(\tpi(\lambda+\rho)).
$$
It remains to compute 
$\iota_F^* {\J}_\lambda(\L^{Car})$. For this we 
we use the description of the Liouville form in the Weil
model. By Proposition \ref{RestrictionKernel} we have 
$$
\iota_F^*\ca{R}_G^T(\L)
=e^{\om_F}\,(\iota_T\circ \Phi_F)^* \ca{R}_G^T(\La)
=e^{\om_F} \Phi_F^* \Lambda_T^-.
$$
Therefore, the image of $\L^{Car}$ under the map 
$\wh{\ca{C}}_G(M)\to \wh{\ca{C}}_T(M)^-$ is given by 
$e^{\om_F}\delta_{\Phi_F}^-$
and we obtain
$$ \iota_F^* \ca{J}_\lambda(\L^{Car})=(\Phi_F)^{\lambda+\rho}\,\,e^{\om_F} 
.$$
\end{proof}

Now we combine the localization formula \eqref{fqdh} and formula
\eqref{interpairings} and use the inverse Fourier transform
on the group $G$ to obtain a formula for intersection pairings 
on reduced spaces,
\beq
\lefteqn{
\int_{M_\Co} \kappa_\Co(\beta) e^{\omega_\Co}} 
\\&&=
k \Vol \mathcal{C} \ \sum_\lambda 
\frac{\chi_\lambda(\mathcal{C}^{-1}) \dim V_\lambda}{(\Vol G)^2}
\ \sum_{F \in \F(\lambda+ \rho)}(\Phi_F)^{\lambda+\rho} 
\int_F  \frac{\iota_F^* \beta_0(2\pi i(\lambda+\rho))}
{\Eul(\nu_F, 2\pi i(\lambda+\rho))}\ 
e^{\omega_F},
\eeq
where $\chi_\lambda(\mathcal{C}^{-1})$ is the value of the character 
$\chi_\lambda$ on the conjugacy class $\mathcal{C}^{-1}$.

\subsection{Generalization to arbitrary compact $G$}
\label{Sec:General}
Up to this point we made the assumption that the compact, connected
group $G$ be a direct product of a simply connected group and a
torus. We will now remove this assumption and show that if
appropriately interpreted, the Duistermaat-Heckman formula
\eqref{fqdh} holds in the general case.

Without this assumption the definition of the Liouville form breaks
down, since the form $\Lambda$ is no longer a well-defined element of
$\wh{W}_G\otimes\Om^*(G)$.  Indeed, in general a Hamiltonian
$G$-manifold with group-valued moment map need not be orientable. 
The simplest counter-example is the conjugacy class in the rotation
group $\SO(3)$ consisting of rotations by an angle $\pi$, which is
isomorphic to $\R P^2$.

Choose a finite covering $\c_G: G' \to G$ such that $G'$ is 
a direct product of a compact, 
connected, simply-connected group and a torus. The kernel
$R:=\c_G^{-1}(e)$ is a subgroup of the center of $G'$. Given a
group-valued Hamiltonian $G$-space $(M, \omega, \Phi)$, consider the
fiber product,
%
\begin{equation*}
\vcenter{\xymatrix{
M' \ar[d]^{\c_M} \ar[r]_{\Phi'} & G' \ar[d]^{\c_G} \\
M \ar[r]_{\Phi} &  G}}
\end{equation*} 
that is
$M' = \{ (x, g') \in M \times G' | \ \Phi(x) =\c_G(g') \}.$
The diagonal action of $G'$ on $M \times G'$ (which in fact is 
a $G$-action since $R$ acts trivially) leaves $M'$ invariant,
and commutes with the action of $R$ by deck transformations. 
Let $\omega'= \c_M^* \omega$. Then $(M', \omega', \Phi')$ 
is a group-valued Hamiltonian $G'$-space. The form $\om'$ is 
invariant with respect to deck transformations, 
while the moment map $\Phi'$ satisfies 
$c^*\Phi'=c\Phi$ for $c\in R$. Given a regular value $g'\in G'$ 
of $\Phi'$, the image $g=\c_G(g')\in G$ is a regular value for 
$\Phi$ and there is a canonical isomorphism $M'_{g'}=M_g$. 

Suppose $\beta\in {\ca{C}}_G(M)$ is an equivariant cocycle, 
and  $\beta'=\c_M^*\beta\in \ca{C}_{G'}(M')$. 
The form $\beta'$ defines a twisted 
DH-distribution $ (\m')^{\beta'}\in
\ca{E}'(G')^{G'}$, invariant under the action of $R$. 
Let a twisted DH-distribution for $M$ be defined by 
$$\m^\beta:=\f{1}{\# R} (\c_M)_*(\m')^{\beta'}\in\ca{E}'(G)^G.$$
It is shown in \cite{al:du} that formula \eqref{interpairings},
relating the values of $\m^\beta$ and intersection pairings on reduced
spaces, carries over to this more general situation.

The Fourier coefficients of $\m^\beta$ can be read off from those 
for $(\m')^{\beta'}$. Given a dominant weight $\lambda\in\Lambda^*_+$,
the corresponding irreducible character $\chi'_\lambda$ for $G'$  
is the pull-back $\c_G^*\chi_\lambda$, hence
$$
\l \mathfrak{m}^\beta, \chi_\lambda \r =
\frac{1}{\# R} \ \l (\mathfrak{m}')^{\beta'}, \chi'_\lambda \r,
$$
The pairing on the right hand side is given by the
localization formula \eqref{fqdh}. For any fixed point manifold
$F\in \F(\lambda+\rho)$ let $F'=\c_M^{-1}(F)$ be its pre-image 
in $M'$. The union of all $F'$ is the fixed point set for 
$(\lambda+\rho)_{M'}$. Therefore, \eqref{fqdh} shows that  
$$
\l \mathfrak{m}^\beta, \chi_\lambda \r = \frac{\dim V_\lambda}{\# R}\
\sum_{F \in \F(\lambda+\rho) }(\Phi')^{\lambda+\rho} \int_{F'}
\frac{\iota_{F'}^*(\beta'(2\pi i (\lambda+\rho))
e^{\omega'} }{\Eul(\nu_{F'}, 2\pi i
(\lambda+\rho))}.
$$
We want to replace integrals over $F'$ by integrals over $F$. Note first
that $(\Phi'_{F'})^{\lambda}$  is just the pull-back of 
$\Phi_F^\lambda$, and that $\iota_{F'}^*\om$ is the pull-back of 
$\iota_F^*\om$. On the other hand, the locally constant 
function $(\Phi'_{F'})^\rho:\,F'\to S^1$ does not 
descend to $F$ since for $c\in R$, 
$$c^*(\Phi'_{F'})^\rho =c^\rho(\Phi'_{F'})^\rho .$$
Since $2\rho\in\Lambda^*$ is a weight for $G$, $c^{2\rho}=1$ so that
$c^\rho=\pm 1$ is just a sign. As shown in Section 8 of 
\cite{al:du}, this sign also gives the change in orientation of 
$M'$ under the action of $c$. It follows that the quotient 
$$
\frac{(\Phi')_{F'}^{\rho} }{\Eul(\nu_{F'}, 2\pi i
(\lambda+\rho))}
$$
is invariant under $R$, and therefore descends to a form on 
$F$. If we denote this form by 
$$
\f{\Phi_F^\rho}{\Eul(\nu_F, 2\pi i (\lambda+\rho))},
$$
then \eqref{fqdh} holds for any compact, connected Lie group.
(Put differently, while $\nu_F$ need not be orientable, 
the Euler form makes sense as a form with values in the orientation
bundle of $\nu_F$, and 
$\Phi_F^\rho$  descends to a section of the orientation bundle.
Their quotient is a form in the usual sense.)

\bibliographystyle{amsplain}
\bibliography{ref}

\end{document}